\documentclass[a4paper, draft]{article}

\usepackage[leqno]{amsmath}
\usepackage{amssymb} %Ez kell $\mathbb R$-hez.
\usepackage{amsthm}  %Ez meg a proof kornyezethez.
\usepackage{hhline}
\usepackage{array}

\newtheorem{tm}{Theorem}[section]
\newtheorem{lm}[tm]{Lemma}
\newtheorem{pr}[tm]{Proposition}
\newtheorem{cor}[tm]{Corollary}
\newtheorem{con}[tm]{Condition}
\newtheorem{rem}[tm]{Remark}

\newcommand{\fej}[1]{\section{\!\!\!\!\!\!.\ #1}}
\newcommand{\alfej}[1]{\subsection{\!\!\!#1}}

\renewcommand{\arraystretch}{1.6} %Tablazatsorvastagitas.

\renewcommand*{\int}{\intop\limits}
\newcommand*{\Zb}{\mathbb Z}
\newcommand*{\Rb}{\mathbb R}
\newcommand*{\om}{\omega}
\newcommand*{\ze}{\zeta}
\newcommand*{\si}{\sigma}
\newcommand*{\te}{\theta}
\newcommand*{\de}{\delta}
\newcommand*{\tom}{\widetilde{\omega}}
\newcommand*{\wh}{\widetilde h}
\newcommand*{\un}[1]{\underline{#1}}
\newcommand*{\al}{\alpha}

\newcommand*{\ga}{\gamma}
\newcommand*{\e}[1]{\text{\rm e}^{#1}}
\newcommand*{\vr}{\varrho}
\newcommand*{\ve}{\varepsilon}
\newcommand*{\tr}{\widetilde{r}}
\newcommand*{\rsl}{r^{S\,\text{left}}}
\newcommand*{\rsr}{r^{S\,\text{right}}}
\newcommand*{\di}{\,\text{\rm d}}
\newcommand*{\Ff}{\,\mathcal F}
\newcommand*{\Vf}{\mathcal V}
\newcommand*{\Bf}{\mathcal B}
\newcommand*{\Ev}{{\bf E}}
\newcommand*{\Vv}{{\text{\bf Var}}}
\newcommand*{\Pv}{{\bf P}}
\newcommand*{\vp}{\varphi}
\newcommand*{\tv}{\mbox{\boldmath$\tau$}}
\newcommand*{\lc}{\lceil}
\newcommand*{\rc}{\rceil}
\newcommand*{\lf}{\lfloor}
\newcommand*{\rf}{\rfloor}
\newcommand*{\ta}{{\lf Vt\rf}}
\newcommand*{\tf}{{\lc Vt\rc}}
\newcommand*{\omin}{\om^{\text{min}}}
\newcommand*{\omax}{\om^{\text{max}}}

\DeclareMathOperator{\sh}{sinh}
\DeclareMathOperator{\ch}{cosh}

\begin{document}
\setlength{\arraycolsep}{.14em}
\title{Growth fluctuations in a class of deposition models}
\author{M\'arton Bal\'azs\\
\medskip
(Institute of Mathematics, Technical University Budapest)}
\date{}
\maketitle

\begin{abstract}
We compute the growth fluctuations in equilibrium of a wide class of deposition models. These models also serve as general frame to several nearest-neighbor particle jump processes, e.g.\ the simple exclusion or the zero range process, where our result turns to current fluctuations of the particles. We use martingale technique and coupling methods to show that, rescaled by time, the variance of the growth as seen by a deterministic moving observer has the form $|V-C|\cdot D$, where $V$ and $C$ is the speed of the observer and the second class particle, respectively, and $D$ is a constant connected to the equilibrium distribution of the model. Our main result is a generalization of Ferrari and Fontes' result for simple exclusion process. Law of large numbers and central limit theorem are also proven. We need some properties of the motion of the second class particle, which are known for simple exclusion and are partly known for zero range processes, and which are proven here for a type of deposition models and also for a type of zero range processes. 
\bigskip
\begin{center}
{\bf R\'esum\'e}
\end{center}

On compute les fluctuations du grandissement dans l'\'etat d'\'equilibre d'une classe vaste des processus de d\'echarge. Ces processus forment aussi bien un cadre pour quelques mod\`eles des bonds voisins des particules, p.\ e.\ le mod\`ele simple exclusion ou zero range, o\`u notre r\'esultats de\-vi\-ennent des r\'esultats sur les fluctuations du flux des particules. On utilise de m\'ethode martingale et des techniques des couplages pour pr\'esenter que le variance du grandissement, regradu\'e par le temps et vu par un observateur qui avance d\'eterminement \`a une vitesse $V$, a la forme $|V-C|\cdot D$, o\`u $C$ est la v\'elocit\'e de la particule de deuxi\`eme classe, et $D$ est une cons\-tante connect\'ee \`a l'\'etat d'\'equilibre du mod\`ele. Notre r\'esultat principal est une g\'en\'eralisation du r\'esultat de Ferrari et Fontes pour le mod\`ele simple exclusion. La loi des grandes nombres et le th\'eor\`eme de la limite centrale sont aussi d\'emontr\'es. Nous avons besoin de quelques propri\'et\'es du mouvement du particule de deuxi\`eme classe, qui sont connues pour simple exclusion et partiellement pour le mod\`ele zero range, et qui sont d\'emontr\'ees ici pour un type des processus de d\'echarge et pour un type des mod\`eles zero range aussi.
\end{abstract}


\fej{Introduction}
Stochastic deposition models can be used to obtain microscopic description of domain growths, e.g.\ a colony of cells or an infected area of plants. The fluctuation of the growth is itself of great interest. Moreover, these models are in close connection to interacting particle systems, where the particle diffusion corresponds to rescaled surface fluctuation. As it is shown below, an additional feature of deposition models is the possibility of handling antiparticles as well as particles in the particle representation of the process. It has been known \cite{se} for the simple exclusion process, that the current fluctuation is in close connection to the motion of the so-called second class particle, and, divided by time, its variance vanishes for an observer moving with the speed of this particle. In this latter case, Pr\"ahofer and Spohn \cite{spohn} suggest this quantity to be in the order of $t^{2/3}$.

In the present note we consider a wide class of one-dimensional deposition models, parameterized by rate functions describing a column's growth depending on the neighboring columns' relative heights. By monotonicity properties of the rate functions, our models are attractive. For a treatment of these models in a hydrodynamical context, without using attractivity, see T\'oth and Valk\'o \cite{hydro}; T\'oth and Werner \cite{balint}. Following Rezakhanlou \cite{rez}, we first show some conditions on the model in order to have product measures as stationary ones for the process. (By stationarity, we mean time-invariance in this paper.) Our description is general enough to include the asymmetric simple exclusion process, some types of the zero range process, and a family of deposition models, which we call bricklayers' models. In this general frame, we compute the growth fluctuations in order $\mathcal O(t)$, hence generalize the result of Ferrari and Fontes \cite{se}. In the computations we couple two processes, which only differ at one site. This is the position of the so-called {\sl defect tracer}, or also called second class particle. We need law of large numbers and a second moment condition for the position of this extra particle. These have been established for simple exclusion \cite{shock}, but, as far as we know, only $\text{L}^1$-convergence is known for most kinds of zero range processes \cite{rez}. We prove $L^n$-convergence with any $n$ for the defect tracer of the totally asymmetric zero range process and for our new bricklayers' models via various coupling techniques.
\alfej{The model}\label{sc:mod}
The class of models described here is a generalization of the so-called misanthrope process. For $-\infty\le\omin\le0$ and $1\le\omax\le\infty$ (possibly infinite valued) integers, we define
\[
I:\,=\left\{z\in\Zb\,:\,\omin-1<z<\omax+1\right\}
\]
and the phase space 
\[
\Omega=\left\{\un\om=(\om_i)_{i\in\Zb}\ :\ \om_i\in I\right\}=I^{\Zb}.
\]
For each pair of neighboring sites $i$ and $i+1$ of $\Zb$, we can imagine a column built of bricks, above the edge $(i,\,i+1)$. The height of this column is denoted by $h_i$. If $\un{\om}(t)\in\Omega$ for a fixed time $t\in\Rb$ then $\om_i(t)=h_{i-1}(t)-h_i(t)\,\in I$ is the negative discrete gradient of the height of the ``wall''. The growth of a column is described by jump processes. A brick can be added:
\[
\left(\om_i,\,\om_{i+1}\right)\longrightarrow\left(\om_i-1,\,\om_{i+1}+1\right)\ \ \ \text{with rate}\ r(\om_i,\,\om_{i+1}).
\]
Conditionally on $\un\om(t)$, these moves are independent. See fig.\,\ref{fig:elso} for some possible instantaneous changes. For small $\ve$, the conditional expectation of the growth of the column between $i$ and $i+1$ in the time interval $[t,\,t+\varepsilon]$ is $r(\om_i(t),\,\om_{i+1}(t))\cdot\varepsilon+\mathfrak{o}(\ve)$.

The rates must satisfy
\[
r(\omin,\,\cdot\,)\equiv r(\,\cdot\,,\,\omax)\equiv0
\]
whenever either $\omin$ or $\omax$ is finite. We assume $r$ to be non-zero in all other cases. We want the dynamics to smoothen our interface, that is why we assume monotonicity in the following way:
\begin{equation}
r(z+1,\,y)\ge r(z,\,y),\qquad\ r(y,\,z+1)\le r(y,\,z)\label{eq:mon}
\end{equation}
for $y,\,z,\,z+1\in I$.
This means that the higher neighbors a column has, the faster it grows. Our model is hence {\sl attractive}.

We are going to use product property of the model's stationary measure. For this reason, similarly to Rezakhanlou \cite{rez}, we assume that for any $x,\,y,\,z\in I$
\begin{equation}
r(x,\,y)+r(y,\,z)+r(z,\,x)=r(x,\,z)+r(z,\,y)+r(y,\,x),\label{eq:stacifelt}
\end{equation}
and for $\omin<x,\,y,\,z<\omax+1$
\begin{equation}
r(x,y-1)\cdot r(y,z-1)\cdot r(z,x-1)=r(x,z-1)\cdot r(z,y-1)\cdot r(y,x-1).\label{eq:ffelt}
\end{equation}
These two conditions imply product structure of the stationary measure, see section \ref{sc:gibbs}. Equation \eqref{eq:ffelt} is equivalent to the condition $r(y,\,z)=s(y,\,z+1)\cdot f(y)$ for some function $f$ and a symmetric function $s$. 

At time $t$, the interface mentioned above is described by $\un{\om}(t)$. Let $\varphi\,:\,\Omega\to\Rb$ be a finite cylinder function i.e.\ $\varphi$ depends on a finite number of values of $\om_i$. The growth of this interface is a Markov process, with the formal infinitesimal generator $L$:
\begin{equation}
(L\varphi)(\un\om)=\sum_{i\in\Zb}r(\om_i,\,\om_{i+1})\cdot\left[\varphi(\dots,\,\om_i-1,\,\om_{i+1}+1,\,\dots)-\varphi(\un\om)\right].\label{eq:gen}
\end{equation}

When constructing the process rigorously, problems may arise due to the unbounded growth rates. The system being one-component and attractive, we assume that, with appropriate growth conditions on the rates, existence of dynamics on a set of tempered configurations $\widetilde\Omega$ (i.e.\ configurations obeying some restrictive growth conditions) can be established by applying methods initiated by Liggett and Andjel \cite{exi} \cite{and}. Technically we assume that $\widetilde\Omega$ is of full measure w.r.t.\ the canonical Gibbs measures defined in section \ref{sc:gibbs}. In fact this has been proved for some kinds of these models, see below. We do not deal with questions of existence of dynamics in the present paper.

\alfej{Examples}

There are three essentially different cases of these models, all of them are of nearest neighbor type.

\begin{enumerate}
\item {\bf Generalized exclusion processes} are described by our models in case both $\omin$ and $\omax$ are finite.
\begin{itemize}
\item {\bf The totally asymmetric simple exclusion process (SE)} introduced by F.\ Spitzer \cite{spi} is described this way by $\omin=0,\ \omax=1$,
\[
r(\om_i,\,\om_{i+1})=\om_i\cdot(1-\om_{i+1}).
\]
Here $\om_i$ is the occupation number for the site $i$, and $r(\om_i,\,\om_{i+1})$ is the rate for a particle to jump from site $i$ to $i+1$. Conditions (\ref{eq:mon}), \eqref{eq:stacifelt} and \eqref{eq:ffelt} for these rates are satisfied.

\item {\bf A particle-antiparticle exclusion process} is also shown to de\-mon\-stra\-te the generality of the frame described above. Let $\omin=-1,\ \omax=1$. Fix $c$ ({\sl creation}), $a$ ({\sl annihilation}) positive rates with $c\le a/2$. Put
\[
r(0,\,0)=c,\ \ r(0,\,-1)=\frac{a}{2},\ \ r(1,\,0)=\frac{a}{2},\ \ r(1,\,-1)=a,
\]
and all other rates are zero. If $\om_i$ is the number of particles at site $i$, with $\om_i=-1$ meaning the presence of an antiparticle, then this model describes a totally asymmetric exclusion process of particles and antiparticles with annihilation and particle-antiparticle pair creation. These rates also satisfy our conditions.
\end{itemize}
Other generalizations are possible allowing a bounded number of particles (or antiparticles) to jump to the same site. By the bounded jump rates and by nearest-neighbor type of interaction, the construction of dynamics of these processes is well understood, see e.g.\ Liggett \cite{ips}.
\item {\bf Generalized misanthrope processes} are obtained by choosing $\omin>-\infty,\ \omax=\infty$.
\begin{itemize}
\item {\bf The zero range process (ZR)} is included by $\omin=0,\ \omax=\infty$,
\[
r(z,\,y)=f(z)
\]
with an arbitrary $f\,:\,\Zb^+\to\Rb^+$ nondecreasing function and $f(0)=0$. Here $\om_i$ represents the number of particles at site $i$. These rates trivially satisfy conditions \eqref{eq:mon}, \eqref{eq:stacifelt}, \eqref{eq:ffelt}. The dynamics of this process is constructed by Andjel \cite{and} under the condition that the rate function $f$ obeys the growth condition $|f(z+1)-f(z)|\le K$ for some $K>0$ and all $z\ge0$.
\end{itemize}
\item {\bf General deposition processes} are the type of these models where $\omin=-\infty$ and $\omax=\infty$. In this case, the height difference between columns next to each other can be arbitrary in $\Zb$. Hence the presence of antiparticles can not be avoided when trying to give a particle representation of the process.
\begin{itemize}
\item {\bf Bricklayers' models (BL).} Let 
\[
r(z,\,y):\,=f(z)+f(-y)
\]
with the property
\[
f(z)\cdot f(-z+1)=1
\]
for the nondecreasing function $f$ and for any $z\in\Zb$. This process can be represented by bricklayers standing at each site $i$, laying a brick on the column on their left with rate $f(-\om_i)$ and laying a brick to their right with rate $f(\om_i)$. This interpretation gives reason to call these models bricklayers' model. Conditions (\ref{eq:mon}), (\ref{eq:stacifelt}) and \eqref{eq:ffelt} hold for $r$. Similarly to the ZR process, this model is constructed by Booth and Quant \cite{lorna} only in case $|f(z+1)-f(z)|$ is bounded in $\Zb$.
\end{itemize}
\end{enumerate}

\alfej{Translation invariant stationary product measures}\label{sc:gibbs}

We are interested in translation invariant stationary measures for these processes, i.e.\ canonical Gibbs-measures. We construct such measures similarly to Rezakhanlou \cite{rez} of the following form. Fix $f(1)>0$ and define
\begin{equation}
f(z):\,=\frac{r(z,\,0)}{r(1,\,z-1)}\cdot f(1)\label{eq:fdef}
\end{equation}
for $\omin<z<\omax+1$. Then $f$ is a nondecreasing strictly positive function. For $I\ni z>0$ we define 
\[
f(z)!:\,=\prod_{y=1}^zf(y),
\]
while for $I\ni z<0$ let
\[
f(z)!:\,=\frac{1}{\prod\limits_{y=z+1}^0f(y)},
\]
finally $f(0)!:\,=1$. Then we have 
\[
f(z)!\cdot f(z+1)=f(z+1)!
\]
for all $z\in I$. Let
\[
\bar\te:\,=\left\{\begin{array}{ll}\log\left(\liminf\limits_{z\to\infty}\left(f(z)!\right)^{1/z}\right)=\lim\limits_{z\to\infty}\log(f(z))\ \ &,\ \text{if}\ \omax=\infty\\\infty\ \ &,\ \text{else}\end{array}\right.
\]
and
\[
\un\te:\,=\left\{\begin{array}{ll}\log\left(\limsup\limits_{z\to\infty}\left(f(-z)!\right)^{1/z}\right)=\lim\limits_{z\to\infty}\log(f(-z))\ \ &,\ \text{if}\ \omin=-\infty\\-\infty\ \ &,\ \text{else}.\end{array}\right.
\]
By monotonicity of $f$, we have $\bar\te\ge\un\te$. We assume $\bar\te>\un\te$. With a generic real parameter $\te\in\left(\un\te,\,\bar\te\right)$, we define
\[
Z(\te):\,=\sum_{z\in I}\frac{\e{\te z}}{f(z)!}.
\]
Let the product-measure $\un\mu_\te$ have marginals
\begin{equation}
\mu_\te(z)=\un\mu_\te\left\{\un\om\,:\,\om_i=z\right\}:\,=\frac{1}{Z(\te)}\cdot\frac{\e{\te z}}{f(z)!}\label{eq:om}.
\end{equation}
By definition it has the property
\[
\frac{\mu_\te(z+1)}{\mu_\te(z)}=\frac{\e{\te}}{f(z+1)}
\]
which implies
\begin{equation}
r(z+1,\,y-1)\cdot\frac{\mu_\te(z+1)\,\mu_\te(y-1)}{\mu_\te(z)\,\mu_\te(y)}=r(y,\,z)\label{eq:ctszimm}
\end{equation}
due to \eqref{eq:fdef} and \eqref{eq:ffelt}. Hence stationarity of $\un\mu_\te$ follows via \eqref{eq:stacifelt}.

As can be verified, the expectation value $\vr(\te):\,=\Ev_\te(\om_i)$ is a strictly increasing function of $\te$. We introduce its inverse $\te(\vr)$ and the function
\begin{equation}
\mathcal H(\vr):\,=\Ev_{\te(\vr)}\left\{r(\om_i,\,\om_{i+1})\right\},\label{eq:haro}
\end{equation}
playing an important role in hydrodynamical considerations. For the SE model, the construction leads to the well-known Bernoulli pro\-duct-me\-a\-su\-re with mar\-gi\-nals
\[
\begin{array}{rcl}
\mu(1)=&\un\mu\{\un\om\,:\,\om_i=1\}&:\,=\varrho,\\
\mu(0)=&\un\mu\{\un\om\,:\,\om_i=0\}&:\,=1-\varrho
\end{array}
\]
with a real number $\vr$ between zero and one (the density of the particles). In our notations, $-\vr$ describes the average slope of the interface.

For the particle-antiparticle exclusion process, the relative probability of having a particle or an antiparticle as a function of the rates goes as $\sqrt{c/a}$, independently for the sites. The density of particles relative to antiparticles can be set by an arbitrary parameter.

Both for the ZR process and for BL models, it turns out that $f$ defined in \eqref{eq:fdef} and $f$ in the definition of the rates agree. 

It is not hard to show ergodicity of these models, which also implies extremality of the invariant measures $\un\mu_\te$:
\begin{pr}
The processes given in subsection \ref{sc:mod}, distributed according to their stationary measures $\un\mu_\te$ \eqref{eq:om}, are ergodic.
\end{pr}
\begin{proof}
We need to show that any (time-) stationary bounded measurable function defined on the trajectories of the process is constant a.s. By proposition V.2.4 of Neveu \cite{neveu}, this follows once we see that any bounded function $\vp$ on $\widetilde\Omega$ satisfying $P\vp=\vp$ is constant for $\un\mu$-almost all $\un\om$ with the Markov-transition operator $P$. Hence ergodicity of the process follows if $L\vp=0$ implies $\vp(\un\om)=\text{constant}$ for almost all $\un\om\in\widetilde\Omega$. We compute the Dirichlet-form 
\begin{multline*}
-\Ev_\te(\vp\cdot L\vp)=\\
=\frac12\Ev_\te\left\{\sum_{i\in\Zb}r(\om_i,\,\om_{i+1})\cdot\left[\vp(\dots,\,\om_i-1,\,\om_{i+1}+1,\,\dots)-\vp(\un\om)\right]^2\right\}.
\end{multline*}
By positivity of the rates, this shows that assuming $L\,\vp\equiv0$ results in
\[
\vp(\dots,\,\om_i-1,\,\om_{i+1}+1,\,\dots)=\vp(\un\om)
\]
for almost all $\un\om\in\widetilde\Omega$. Consecutive use of this equation shows that any function obeying $L\vp=0$ does a.s.\ not depend on any finite cylinder set in $\widetilde\Omega$. Especially, for $\ve>0$ and a constant $K\in\text{Ran}(\vp)$, the event
\[
\{\vp(\un\om)\in(K,\,K+\ve]\}
\]
does not depend on any finite cylinder set. Hence by Kolmogorov's 0-1 law, the probability of these events is zero or one w.r.t.\ the product measure $\un\mu$. Partitioning the bounded image of $\vp$, this shows that this function is constant for almost all $\un\om$. 
\end{proof}

\alfej{Results}

We start our model in a canonical Gibbs-distribution, with parameter $\te$. For a fixed speed value $V>0$ we define
\[
J^{(V)}(t):=h_\ta(t)-h_0(0),
\]
the height of column at site $\ta$ at time $t$, relative to the initial height of the column at the origin. For $V<0$, we introduce
\[
J^{(V)}(t):=h_\tf(t)-h_0(0),
\]
which is the mirror-symmetric form of $J^{(V)}$ defined above for positive $V$'s. For $V=0$ we write
\[
J(t)=J^{(0)}(t):\,=h_0(t)-h_0(0).
\]
In particle notations of the models, $J^{(V)}(t)$ is the current, i.e.\ the algebraic number of particles jumping through the moving window positioned at $Vt$, in the time interval $[0,\,t]$. We prove law of large numbers for this quantity:
\begin{equation}
\lim_{t\to\infty}\frac{J^{(V)}}{t}=\Ev(r)-V\,\Ev(\om)\ \ \text{a.s.}\label{eq:jfnsztv}
\end{equation}

We need law of large numbers and a second-moment condition for the position $Q(t)$ of the defect tracer (also called second class particle, see section \ref{sc:coupling} for its definition) if one of the coupled models is started from its canonical Gibbs-measure:
\begin{con}\label{con:con}
With initial distribution $\un\mu_\te$ of $\un\om$, weak law of large numbers
\begin{equation}
\lim_{t\to\infty}\Pv_\te\left(\left|\frac{Q(t)}{t}-C(\te)\right|>\delta\right)=0\label{eq:nsztv}
\end{equation}
for a speed value $C(\te)$ and for any $\delta>0$ holds, and the bound
\begin{equation}
\Ev_\te\left(\frac{Q(t)^2}{t^2}\right)<K<\infty\label{eq:masodq}
\end{equation}
is satisfied for all large $t$ for the position $Q(t)$ of the defect tracer.
\end{con}
\noindent
Inequality \eqref{eq:masodq} is obvious in case of bounded rates, since in this situation, the process $|Q(t)|$ is bounded by some Poisson-process. 
\begin{tm}[Main]\label{tm:main} Assume condition \ref{con:con}. Then
\begin{equation}
\lim_{t\to\infty}\frac{\Vv_\te(J^{(V)}(t))}{t}=|V-C(\te)|\cdot\Vv_\te(\om_0)=\,:D_J(\te)\label{eq:fovegso}
\end{equation}
for any $V\in\Rb$, where $\Vv_\te$ stands for the variance w.r.t.\ $\mu_\te$.
\end{tm}
\begin{tm}[Central limit theorem]\label{tm:cht} Assuming condition \ref{con:con},
\[
\lim_{t\to\infty}\Pv_\te\left(\frac{\widetilde J^{(V)}(t)}{\sqrt{D_J(\te)}\cdot\sqrt{t}}\le x\right)=\Phi(x)=\int_{-\infty}^x\frac{\e{-y^2/2}}{\sqrt{2\pi}}\di y,
\]
i.e.\ $\widetilde J^{(V)}(t)/\sqrt{t}$ converges in distribution to $N(0,\,D_J(\te))$, a centered normal random variable with variance $D_J(\te)$ of \eqref{eq:fovegso}. Tilde means here that the mean value of $J^{(V)}(t)$ is subtracted. 
\end{tm}

For the SE model, \eqref{eq:nsztv} is proven in \cite{shock}. It is shown there that 
\[
\lim_{t\to\infty}\frac{Q(t)}{t}=1-2\vr\ \ \ \text{a.s.}
\]
Condition \ref{con:con} is satisfied by this law, hence theorem \ref{tm:main} gives
\[
\lim_{t\to\infty}\frac{\Vv(J^{(V)}(t))}{t}=\vr\,(1-\vr)\,|(1-2\,\vr)-V|,
\]
and the central limit theorem \ref{tm:cht} also holds. These results have been known for SE by Ferrari and Fontes \cite{se}. 

For the ZR and BL models, we need a condition on the growth rates:
\begin{con}\label{con:cvx}
For ZR and BL processes defined above, the rate function $f$ is convex.
\end{con}
\noindent
For the ZR process, under this condition and assuming either strict convexity or concavity of $\mathcal H(\vr)$ defined in \eqref{eq:haro}, more than \eqref{eq:nsztv}, namely, $\text{L}^1$-con\-ver\-gen\-ce is established by Rezakhanlou \cite{rez} with speed
\begin{equation}
C(\te)=\frac{\e\te}{\Vv_\te(\om)}.\label{eq:zrnsztv}
\end{equation}
As far as we know, the second-moment condition \eqref{eq:masodq} has not yet been proven for this model. 
\begin{tm}\label{tm:con}
For ZR and BL models satisfying condition \ref{con:cvx} with initial distribution $\un\mu_\te$ of $\un\om$, and for any $n\in\Zb^+$,
\[
\frac{Q(t)}{t}\to C(\te)\qquad\text{in}\ L^n,
\]
where $C(\te)$ is defined in \eqref{eq:zrnsztv} for the ZR process, and
\begin{equation}
C(\te):\,=\frac{2\sh(\te)}{\Vv_\te(\om)}\label{eq:seb}
\end{equation}
for the BL model.
\end{tm}
Hence under condition \ref{con:cvx}, condition \ref{con:con} and thus theorem \ref{tm:main} and \ref{tm:cht} hold for both ZR and BL models with $C(\te)$ defined in \eqref{eq:zrnsztv} and \eqref{eq:seb}, respectively.
As we expect by mirror symmetric properties of the BL model, the speed $C(\te)$ of the defect tracer is zero in case $\te=0$ in this model. 

Our methods do not rely on hydrodynamic limits. $C(\te)$ is a nondecreasing function for the totally asymmetric ZR process and BL model under condition \ref{con:cvx}, see remark \ref{rm:dupla}. This shows (non strict) convexity of the function $\mathcal H(\vr)$ of \eqref{eq:haro} for these models, since 
\[
C(\te(\vr))=\frac{\di\mathcal H(\vr)}{\di\vr}
\]
after some computations, and $\te(\vr)$ is also a monotone function. 
\begin{pr}\label{pr:cvx}
Under condition \ref{con:cvx}, the function $\mathcal H(\vr)$ is strictly convex for the BL model. For the ZR process satisfying \ref{con:cvx}, linearity of $\mathcal H(\vr)$ is equivalent to linearity of the rate function $f$ on $\Zb$, which is the case of independent random walk of the particles. If this is not the case, then $\mathcal H(\vr)$ is strictly convex.
\end{pr}
This is an important observation for \cite{valak}, since this property is only proved for small $\te$ values there. It is also remarkable for \cite{rez}, where strict convexity is just assumed.

We remark that rates for removal of the bricks can also be introduced to obtain a model with both growth and decrease of columns. In particle notations this represents possible left jumps of particles (or right jump of antiparticles, respectively). Therefore, not only the totally asymmetric case, but the general asymmetric case of particle processes (SE or ZR, for example) can also be included in the description. The extension of the proof of theorems \ref{tm:main} and \ref{tm:cht} to this case is straightforward. However, the coupling arguments used to establish condition \ref{con:con} for ZR and BL models in later sections are not applicable in case of brick-removal.

We see that $\lim\limits_{t\to\infty}\frac{\Vv(J^{(V)}(t))}{t}$ vanishes if we observe this quantity from the moving position $Vt=C(\te)t$, having the characteristic speed of the hydrodynamical equation. This has been known for the SE model with strongly restricted values of $\om_i$, and now it is proven for the class of more general models with possibly $\om_i\in\text{all}\,\Zb$ also. The interesting question, of which the answer is strongly suggested for some models \cite{spohn}, is the correct exponent of $t$ leading to nontrivial limit of $\Vv(J^{(C)}(t))/t^{2\al}$ as $t\to\infty$. $\al$ is believed to be 1/3, in close connection to $t^{2/3}$ order fluctuations of the position $Q(t)$ of the defect tracer.

The structure of the paper is the following: after some definitions on the reversed chain, we begin with separating martingales from $\Vv(J(t))$ in section \ref{sc:mart}. Then we proceed in section \ref{sc:inv} by computing the generator's inverse on the rates and then by transforming $\Vv(J(t))$ into nontrivial correlations. These correlations can be computed using monotonicity thus coupling possibilities of the model, this is done in section \ref{sc:coupling}. This section also includes a technical lemma showing an interesting relation of space-time correlations to the motion of the defect tracer. After $J(t)$, we deal with $J^{(V)}(t)$, the growth in non-vertical directions in section \ref{sc:sl}. Our results are proven in this section, except for theorem \ref{tm:con}, which is proven in the last section for the totally asymmetric ZR process and for BL models. This last section includes the introduction of a new random walk depending on our processes, and new coupling techniques based on convexity of the rate function $f$. As another consequence of these methods, this part is followed by a proof of strict convexity of the function $\mathcal H(\vr)$.

\fej{The growth and correlations}

In this section we obtain a formula for $\Vv(J(t))$, which contains only space-time correlations of $\om_i(t)$'s as non-trivial expressions. 

\alfej{The reversed chain}

The formal infinitesimal generator $L^*$ for the reversed chain is of the form
\[
(L^*\varphi)(\un\om)=\sum_{i\in\Zb}r^*(\om_i,\,\om_{i+1})\cdot\left[\varphi(\dots,\,\om_i+1,\,\om_{i+1}-1,\,\dots)-\varphi(\un\om)\right]\label{eq:lcsillag}
\]
on the finite cylinder functions. The rates $r^*$ of the reversed process w.r.t.\ $\un\mu_\te$ can be determined by the equation
\[
\Ev_\te\left(\psi(\un\om)\cdot L\vp(\un\om)\right)=\Ev_\te\left(\vp(\un\om)\cdot L^*\psi(\un\om)\right).
\]
\begin{pr}For $\omin\le z,\,y\le\omax$,
\begin{equation}
r^*(z,\,y)=r(y,\,z).\label{eq:revr}
\end{equation}
\end{pr}
Note that the rates of the reversed process do not depend on the parameter $\te$ of the original process' distribution.
\begin{proof}
Let $\psi,\ \vp$ be finite cylinder functions, and let $\mathcal I\subset\Zb$ be a finite discrete interval of which the size can be divided by three, and which contains the set
\[
\left\{i\in\Zb\,:\,\psi,\ \text{or}\ \vp\ \text{depends on}\ \om_i\ \text{or on}\ \om_{i-1}\right\}.
\]
Then the summation index $i$ in the definition \eqref{eq:gen} of the generator can be run on the set $\mathcal I$. We begin by changing variables $\om_i,\,\om_{i+1}$:
\begin{multline*}
\Ev_\te\left(\psi(\un\om)\cdot L\vp(\un\om)\right)=\\
=\Ev_\te\sum_{i\in\mathcal I}\left\{r(\om_i,\,\om_{i+1})\cdot\left[\psi(\un\om)\vp(\dots,\,\om_i-1,\,\om_{i+1}+1,\,\dots)-\psi(\un\om)\cdot\vp(\un\om)\right]\right\}=\\
=\Ev_\te\sum_{i\in\mathcal I}\biggl\{r(\om_i+1,\,\om_{i+1}-1)\cdot\frac{\mu_\te(\om_i+1)\,\mu_\te(\om_{i+1}-1)}{\mu_\te(\om_i)\,\mu\te(\om_{i+1})}\times\\
\times\psi(\dots,\,\om_i+1,\,\om_{i+1}-1,\,\dots)\vp(\un\om)\biggr\}-\Ev_\te\left\{\left(\sum_{i\in\mathcal I}r(\om_i,\,\om_{i+1})\right)\cdot\psi(\un\om)\vp(\un\om)\right\}.
\end{multline*}
Since $|\mathcal I|$ can be divided by three, we can apply \eqref{eq:stacifelt} in order to show that
\[
\sum_{i\in\mathcal I}r(\om_i,\,\om_{i+1})=\sum_{i\in\mathcal I}r(\om_{i+1},\,\om_i)
\]
in the second term. By using \eqref{eq:ctszimm} for the first term we finally obtain
\begin{multline*}
\Ev_\te\left(\psi(\un\om)\cdot L\vp(\un\om)\right)=\\
=\Ev_\te\sum_{i\in\mathcal I}\left\{r(\om_{i+1},\,\om_i)\cdot\left[\psi(\dots,\,\om_i+1,\,\om_{i+1}-1,\,\dots)\vp(\un\om)-\psi(\un\om)\cdot\vp(\un\om)\right]\right\},
\end{multline*}
which equals to $\Ev_\te\left(\vp(\un\om)\cdot L^*\psi(\un\om)\right)$ by choosing $r^*$ according to \eqref{eq:revr}.
\end{proof}
Combining \eqref{eq:ctszimm} with \eqref{eq:revr} leads to
\begin{equation}
r^*\left(z,\,y\right)=\frac{\mu_\te(z+1)\,\mu_\te(y-1)}{\mu_\te(z)\,\mu_\te(y)}\cdot r\left(z+1,\,y-1\right),\label{eq:rrev}
\end{equation}
which is the natural formula suggested by considering conditional expectation values.

In order to simplify notations, let
\[
r(t):\,=r(\om_0(t),\,\om_1(t)),\qquad\qquad r^*(t):\,=r^*(\om_0(t),\,\om_1(t)).
\]

\alfej{Preparatory computations}\label{sc:mart}

For a quantity $A(\un\om)$ with $\Ev|A|<\infty$, let
\[
\widetilde{A}=\widetilde{A}(\un\om):\,=A-\Ev A.
\]
\begin{lm}\label{lm:mart}
$\Vv(J(t))=t\,\Ev(r)+2\int_0^t(t-v)\,\Ev\left(\tr(v)\,r^*(0)\right)\di v.$
\end{lm}
\begin{proof}
By definition, $\Ev(J(t)\,|\,\un\om(0))=t\,r(0)+\mathfrak{o}(t)$, hence
\[
M(t):\,=J(t)-\int_0^tr(s)\di s
\]
is a martingale with $M(0)=0$. Using this,
\begin{equation}
\Vv(J(t))=\Ev M(t)^2+2\,\Ev\left(M(t)\,\int_0^t\tr(s)\di s\right)+\Ev\left(\left(\int_0^t\tr(s)\di s\right)^2\right).\label{eq:rovidebb}
\end{equation}
Due to $\Ev\left(M(t)^2\,|\,\un\om(0)\right)=t\,r(0)+\mathfrak{o}(t)$, the process
\[
N(t):\,=M(t)^2-\int_0^t r(s)\di s
\]
is also a martingale with $N(0)=0$. Hence
\[
\Ev M(t)^2=t\,\Ev(r).
\]
Using the martingale property of $M$, the second term of \eqref{eq:rovidebb} can be written as
\[
2\,\int_0^t\Ev\left(M(t)\,\tr(s)\right)\di s=2\int_0^t\Ev\left(M(s)\,\tr(s)\right)\di s.
\]
Simply changing the limits of integration in the third term of \eqref{eq:rovidebb}, we have
\[
\Ev\left(\left(\int_0^t\tr(s)\di s\right)^2\right)=2\int_0^t\Ev\left(\tr(s)\int_0^s\tr(u)\di u\right)\di s.
\]
These calculations lead to
\begin{multline}
\Vv(J(t))=t\,\Ev(r)+2\int_0^t\Ev\left(\tr(s)\left(M(s)+\int_0^s\tr(u)\di u\right)\right)\di s=\\
=t\,\Ev(r)+2\int_0^t\Ev\left(\tr(s)\,J(s)\right)\di s.\label{eq:hosszu}
\end{multline}

In order to handle $\Ev\left(\tr(s)\,J(s)\right)$, we introduce $J^{(s)\,*}$, the quantity corresponding to $J$ in the reversed model by
\[
J^{(s)\,*}(u):\,=J(s)-J(s-u)\ \ \ \ \ (s\ge u\ge0). 
\]
This is the number of bricks removed from the column in the reversed model started from time $s$. As in case of $J(t)$, a reversed martingale can be separated by
\[
M^{(s)\,*}(u):\,=J^{(s)\,*}(u)-\int_0^ur^*(s-v)\di v.
\]
For this reversed object,\,\ $M^{(s)\,*}(0)=0$\,\ and\,\ $\Ev\left(M^{(s)\,*}(u)\,|\,\Ff_{[t,\,\infty)}\right)=M^{(s)\,*}(s-t)$\ \ if\ \ $0\le s-t\le u$, where $\Ff$ stands for the natural filtration of the (forward) process. In view of this, 
\begin{multline*}
\Ev\left(\tr(s)\,J(s)\right)=\Ev\left[\tr(s)\,\Ev\left(J^{(s)\,*}(s)\,|\,\Ff_{[s,\,\infty)}\right)\right]=\\
=\Ev\left(\tr(s)\int_0^sr^*(s-v)\di v\right)=\int_0^s\Ev\left(\tr(v)\,r^*(0)\right)\di v,
\end{multline*}
where in the last step we used time-invariance of the measure. Using this result, we obtain
\[
\Vv(J(t))=t\,\Ev(r)+2\int_0^t(t-v)\ \Ev\left(\tr(v)\,r^*(0)\right)\di v
\]
from (\ref{eq:hosszu}) by changing the order of integration.
\end{proof}

\alfej{Occurrence of space-time correlations}\label{sc:inv}

In this subsection we denote $r(\om_i,\,\om_{i+1})$ and $\tr(\om_i,\,\om_{i+1})$ by $r_i$ and $\tr_i$, respectively.
For $k\in\Zb$, let 
\[
d_k\ \ :\ \ \Omega\to I\ \ ;\ \ d_k(\un\om)=\om_k
\]
be the $k$-th coordinate of $\Omega$. Then
\begin{equation}
\begin{array}{rll}
\left(Ld_k\right)(\un\om)&=r_{k-1}-r_k&\text{and}\\
\left(L^*d_k\right)(\un\om)&=-r^*_{k-1}+r^*_k&,
\end{array}\label{eq:szomsz}
\end{equation}
where $L^*$ is the infinitesimal generator (\ref{eq:lcsillag}) for the reversed process.
\begin{lm}\label{lm:elketto}
For $0<\al<1$ the expressions
\begin{equation}
\vp_\al:\,=\sum_{k=1}^{\infty}\al^{k-1}d_k\qquad\qquad\qquad\psi_\al:\,=\sum_{k=0}^{\infty}\al^kd_{-k}\label{eq:fidef}
\end{equation}
exist a.s., and
\begin{alignat*}{4}
\lim_{\al\to1}(L\vp_\al)(\un\om)&=-\lim_{\al\to1}(L\psi_\al)(\un\om)&=&\,\tr,\\
\lim_{\al\to1}(L^*\psi_\al)(\un\om)&=-\lim_{\al\to1}(L^*\vp_\al)(\un\om)&=&\,\widetilde{r^*}
\end{alignat*}
in $\text{L}^2$.
\end{lm}
\begin{proof}
The a.s.\ existence of the sums above can be easily shown by using the Borel-Cantelli lemma for the sets 
\[
A_n:\,=\left\{\un\om\ :\ |\om_n|\ge n\right\}.
\]
We show the first equation for $\vp_\al$. By (\ref{eq:szomsz}) 
\begin{multline}
\left(L\vp_\al\right)(\un\om)=r_0+(\al-1)\,\sum_{k=1}^{\infty}r_k\al^{k-1}=\\
=r_0-\Ev\,r+(\al-1)\,\sum_{k=1}^{\infty}(r_k-\Ev\,r)\al^{k-1}=\tr_0+(\al-1)\,\sum_{k=1}^{\infty}\tr_k\al^{k-1}.\label{eq:l2}
\end{multline}
By independence of $\om_i$ and $\om_j$ for $i\ne j,\ \Ev(\tr_l\cdot\tr_k)=0$ if $|l-k|>1$ and $\left|\Ev(\tr_l\cdot\tr_k)\right|\le\Ev(\tr_l\cdot\tr_l)=||\tr||_2^2$, if $|k-l|=0$ or 1. Hence the $\text{L}^2$-norm of the second term on the right-hand side of (\ref{eq:l2}) tends to zero as $\al\to1$:
\begin{multline*}
\left|\left|(\al-1)\,\sum_{k=1}^{\infty}\tr_k\al^{k-1}\right|\right|_2^2
\le(\al-1)^2\sum_{k=1}^{\infty}||\tr_k||_2^2\al^{2k-2}+\\
+2(\al-1)^2\sum_{k=1}^{\infty}||\tr_k||_2^2\al^{2k-3}
=\frac{(\al-1)^2}{1-\al^2}||\tr||_2^2(1+2\al^{-1})\underset{\al\to1}{\longrightarrow}0.
\end{multline*}
The proof of the other three equations is similar.
\end{proof}
Now we can compute the integrals in our expression for $\Vv(J)$.
\begin{tm}\label{tm:formula1}
\begin{eqnarray}
\Vv(J(t))&=&t\,\Ev(r)-2\,t\,\Ev(r^*(0)\cdot\widetilde\om_1(0))+2\sum_{n=1}^{\infty}n\Ev(\tom_0(0)\,\tom_n(t))=\nonumber\\
&=&t\,\Ev(r)+2\,t\,\Ev(r^*(0)\cdot\widetilde\om_0(0))+2\sum_{n=1}^{\infty}n\Ev(\tom_0(0)\,\tom_{-n}(t)).\label{eq:formula1}
\end{eqnarray}
\end{tm}
As can be seen in the next session, the sums on the right-hand side are convergent.
\begin{proof}
Using $\text{L}^2$ convergence stated in lemma \ref{lm:elketto} and Cauchy's inequality, we rewrite the integral in the result of lemma \ref{lm:mart}. We can write $\tr^*(0)$ instead of $r^*(0)$ there, since $\Ev(\widetilde A\,B)=\Ev(\widetilde A\,\widetilde B)$ if both sides exist. 
\begin{multline*}
\left|\int_0^t(t-v)\Ev(\tr(v)\,\tr^*(0))\di v-\lim_{\al\to1}\int_0^t(t-v)\Ev(L\vp_\al(v)\,\tr^*(0))\di v\right|\le\\
\le\lim_{\al\to1}\int_0^t(t-v)\sqrt{\Ev\left([\tr(v)-L\vp_\al(v)]^2\right)\cdot\Ev(\tr^*(0)^2)}\di v=\\
=\lim_{\al\to1}\sqrt{\Ev\left([\tr(0)-L\vp_\al(0)]^2\right)\cdot\Ev(\tr^*(0)^2)}\,\int_0^t(t-v)\di v=0,
\end{multline*}
hence we can apply integration by parts:
\begin{multline*}
\Vv(J(t))=t\,\Ev(r)+2\lim_{\al\to1}\int_0^t(t-v)\,\Ev\left(L\vp_\al(v)\,\tr^*(0)\right)\di v=\\
=t\,\Ev(r)+2\lim_{\al\to1}\int_0^t(t-v)\,\frac{\di}{\di v}\Ev(\vp_\al(v)\,\tr^*(0))\di v=\\
=t\,\Ev(r)-2\,t\lim_{\al\to1}\Ev(\widetilde\vp_\al(0)\,\tr^*(0))+2\lim_{\al\to1}\int_0^t\Ev(\widetilde\vp_\al(v)\,\tr^*(0))\di v.
\end{multline*}
The last integral here can be transformed in the same way, using lemma \ref{lm:elketto} again:
\begin{multline*}
\int_0^t\Ev(\widetilde\vp_\al(v)\,\tr^*(0))\di v=\int_0^t\Ev(\widetilde\vp_\al(0)\,\tr^*(-v))\di v=\\=
\lim_{\ga\to1}\int_0^t\Ev(\widetilde\vp_\al(0)\,L^*\psi_\ga(-v))\di v=\lim_{\ga\to1}\int_0^t\frac{\di}{\di v}\Ev(\widetilde\vp_\al(0)\,\psi_\ga(-v))\di v=\\
=\lim_{\ga\to1}\Ev(\widetilde\vp_\al(0)\,\widetilde\psi_\ga(-t))-\lim_{\ga\to1}\Ev(\widetilde\vp_\al(0)\,\widetilde\psi_\ga(0)).
\end{multline*}
Hence with definitions \eqref{eq:fidef}, the variance of $J(t)$ can now be written as
\begin{multline*}
\Vv(J(t))=t\,\Ev(r)-2\,t\lim_{\al\to1}\Ev(\widetilde\vp_\al(0)\,\tr^*(0))+\\
+2\lim_{\al,\ga\to1}\Ev(\widetilde\vp_\al(0)\,\widetilde\psi_\ga(-t))-2\lim_{\al,\ga\to1}\Ev(\widetilde\vp_\al(0)\,\widetilde\psi_\ga(0))=\\
=t\,\Ev(r)-2\,t\,\lim_{\al\to1}\Ev\left(\sum_{k=1}^{\infty}\al^{k-1}\,\widetilde\om_k(0)\,\tr^*(0)\right)+\\
+2\,\lim_{\al,\ga\to1}\Ev\left(\sum_{k=1}^{\infty}\al^{k-1}\,\widetilde\om_k(0)\sum_{l=0}^{\infty}\ga^l\,\widetilde\om_{-l}(-t)\right)-\\
-2\,\lim_{\al,\ga\to1}\Ev\left(\sum_{k=1}^{\infty}\al^{k-1}\,\widetilde\om_k(0)\sum_{l=0}^{\infty}\al^l\,\widetilde\om_{-l}(0)\right).
\end{multline*}
Using product property of the measure at time $t=0$ and the fact that $r^*$ depends only on $\om_0$ and $\om_1$, most of our expressions become simple (recall that all quantities with tilde are centered random variables):
\begin{multline*}
\Vv(J(t))
=t\,\Ev(r)-2\,t\,\Ev(\tom_1(0)\,\tr^*(0))+2\sum_{k=1}^\infty\sum_{l=0}^\infty\Ev(\tom_k(0)\,\tom_{-l}(-t))-0=\\
=t\,\Ev(r)-2\,t\,\Ev(\tom_1(0)\,\tr^*(0))+2\sum_{n=1}^\infty n\,\Ev(\tom_n(t)\,\tom_0(0)).
\end{multline*}
In the last step, we used translation- and time-invariance of the measure.

We needed $L\vp_\al\to\tr$ and $L^*\psi_\al\to\tr^*$ in $\text{L}^2$ so far. The properties $-L\psi_\al\to\tr$ and $-L^*\vp_\al\to\tr^*$ can be used in a similar way to prove the second equation of the theorem. However, we need both $\vp_\al$ and $\psi_\al$: using only one of them would have lead to a divergent sum in the last step.
\end{proof}
The first two expressions of formula \eqref{eq:formula1} can be computed easily. The difficulty is in determining the space-time correlations $\Ev(\widetilde\om_0(0)\,\widetilde\om_k(t))$. In order to do this, we use coupling technique.

\fej{Coupling and correlations}\label{sc:coupling}

In this section, we show how to couple a pair of our models, with the help of the so-called second class particles. We can use second particles to compute our expressions containing space-time correlations.

\alfej{The basic coupling}\label{sc:bas}
We consider two realizations of a model, namely, $\un\ze$ and $\un\eta$. We show the basic coupling preserving
\begin{equation}
\ze_i(t)\ge\eta_i(t),\label{eq:rend}
\end{equation}
if this property holds initially for $t=0$. We say that $n=\ze_i(t)-\eta_i(t)\ge 0$ is the number of {\sl second class particles} present at site $i$ at time $t$. During the evolution of the processes, the total number of these particles is preserved, and each of them performs a nearest neighbor random walk.

The height of the column of $\un\ze$ (or $\un\eta$) between sites $i$ and $i+1$ is denoted by $g_i$ (or $h_i$, respectively). (These quantities are just used for easier understanding, they are not essential for the processes.) Let $g_i\uparrow$ (or $h_i\uparrow$) mean that the column of $\un\ze$ (or the column of $\un\eta$, respectively) between the sites $i$ and $i+1$ has grown by one brick. Then the coupling rules are shown in table \ref{tab:bas}. Each line of this table represents a possible move, with rate written in the first column. In the last column, $\curvearrowright$ (or $\curvearrowleft$) means that a second class particle has jumped from $i$ to $i+1$ (or from $i+1$ to $i$, respectively). This coupling for the SE model is described (with particle notations) in Liggett \cite{couse}, \cite{ips} and \cite{stochi}. The rates of these steps are non-negative due to \eqref{eq:rend} and monotonicity \eqref{eq:mon} of $r$. These rules clearly preserve property \eqref{eq:rend}, since the rate of any move which could destroy this condition becomes zero. Summing the rates corresponding to either $g_i\uparrow$ or to $h_i\uparrow$ shows that each $\un\ze$ and $\un\eta$ evolves according to its own rates. It would be possible to couple models possessing rates for removal of bricks as well.

\alfej{Correlations and the defect tracer}

We introduce the notation $\un\de_i\in\Omega$, a configuration being one at site $i$ and zero at all other sites. Let $\un\om$ be a model distributed according to $\un\mu_\te$, and $\un\ze(0)=\un\om(0)+\un\de_0$, i.e.\ we have a single one second class particle between $\un\ze$ and $\un\om$, initially at site 0. In order to avoid confusions, we call this particle the {\sl defect tracer}. According to the basic coupling, this single defect tracer is conserved for any time $t$:
\begin{equation}
\un\ze(t)=\un\om(t)+\un\de_{Q(t)}\label{eq:delt}
\end{equation}
The quantity $Q(t)$ is the position of the defect tracer, performing a nearest neighbor random walk on $\Zb$. 

In this subsection we consider the process $(\un\om(t),\,Q(t))$, the model distributed according to the Gibbs measure $\un\mu$ and the random walk $Q(t)$ connected to it with $Q(0)=0$. Using condition \ref{con:con}, we prove theorem \ref{tm:main} for $V=0$. We begin with a technical lemma, showing how to make use of the defect tracer.
\begin{lm}
For the pair $(\un\om(t),\,Q(t))$ defined above, and for a function $F\,:\,I\to\Rb$ with $F(\omax)=0$ and with finite expectation value $\sum F(z)\,\mu(z)$,
\begin{multline}
\Ev\left(\om_n(t)\,\left[\frac{F(\om_0(0)-1)\,\mu(\om_0(0)-1)}{\mu(\om_0(0))}-F(\om_0(0))\right]\right)=\\
=\Ev\left({\bf1}\{Q(t)=n\}\,F(\om_0(0))\right).\label{eq:prefo}
\end{multline}
\end{lm}
\begin{proof}
We take conditional expectation value of \eqref{eq:delt}:
\begin{equation}
\Ev\left(\ze_n(t)\,|\,\om_0(0)=z\right)=\Ev\left(\om_n(t)\,|\,\om_0(0)=z\right)+\Pv\left(Q(t)=n\,|\,\om_0(0)=z\right).\label{eq:lep}
\end{equation}
Initially, $\un\ze(0)=\un\om(0)+\un\de_0$. Therefore, $\un\ze$ itself is also a model with initial distribution $\un\mu$, except for the origin. Hence
\[
\Ev\left(\ze_n(t)\,|\,\om_0(0)=z\right)=\Ev\left(\ze_n(t)\,|\,\ze_0(0)=z+1\right)=\Ev\left(\om_n(t)\,|\,\om_0(0)=z+1\right),
\]
and \eqref{eq:lep} can be written as
\[
\Ev\left(\om_n(t)\,|\,\om_0(0)=z+1\right)-\Ev\left(\om_n(t)\,|\,\om_0(0)=z\right)=\Pv\left(Q(t)=n\,|\,\om_0(0)=z\right).
\]
We multiply both sides with $F(z)\,\mu(z)$ and then add up for all $z\in I$ to obtain
\begin{multline*}
\sum_{z\in I}\Ev\left(\om_n(t)\,|\,\om_0(0)=z\right)\cdot\left(F(z-1)\,\mu(z-1)-F(z)\,\mu(z)\right)=\\
=\sum_{z\in I}\Pv\left(Q(t)=n\,|\,\om_0(0)=z\right)\cdot F(z)\,\mu(z).
\end{multline*}
Here we used that $F(\omax)=0$ and we write $\mu(\omin-1)=0$. We know that $\Pv(\om_0(0)=z)=\mu(z)$, hence the proof follows.
\end{proof}
\begin{cor}\label{cor:g-s}
We use the convention that the empty sum equals zero. Let
\[
g(z):\,=z-\sum_{y\in I}y\,\mu(y).
\]
For $n\in\Zb$, 
\[
\Ev(\widetilde\om_0(0)\,\widetilde\om_n(t))=\Ev\left({\bf1}\{Q(t)=n\}\cdot\sum_{z=\om_0+1}^{\omax}g(z)\frac{\mu(z)}{\mu(\om_0)}\right).
\]
\end{cor}
\begin{proof}
By the previous lemma, our goal is now to find the correct function $F$, for which 
\[
\frac{F(z-1)\,\mu(z-1)}{\mu(z)}-F(z)=g(z)=z-\sum_{y\in I}y\,\mu(y)
\]
is satisfied. By inverting the operation on the left side, we find
\[
F(z):\,=\sum_{y=z+1}^{\omax}g(y)\,\frac{\mu(y)}{\mu(z)}.
\]
This function satisfies the conditions of the lemma. Using (\ref{eq:prefo}),
\begin{multline*}
\Ev(\widetilde\om_n(t)\,\widetilde\om_0(0))=\Ev(\om_n(t)\,\widetilde\om_0(0))=\Ev(\om_n(t)\cdot g(\om_0(0)))=\\=\Ev\left(\om_n(t)\,\left[\frac{F(\om_0(0)-1)\,\mu(\om_0(0)-1)}{\mu(\om_0(0))}-F(\om_0(0))\right]\right)=\\
=\Ev\left({\bf1}\{Q(t)=n\}\,F(\om_0(0))\right)=\\
=\Ev\left({\bf1}\{Q(t)=n\}\cdot\sum_{y=\om_0(0)+1}^{\omax}g(y)\frac{\mu(y)}{\mu(\om_0(0))}\right).
\end{multline*}
\end{proof}

Now it becomes clear that we need to know something about the motion of the defect tracer. $\un\ze$ and $\un\om$ can not be started together from their original stationary distribution due to the initial difference between them, present at the origin. We could follow our defect tracer. Knowing a measure stationary as seen from site $Q(t)$ for all time $t$ would help us to state the law of large numbers for the $Q(t)$ process. In general, we don't know such a stationary measure which has the same asymptotics far on the left and far on the right side. It is shown in \cite{valak}, that under some weak assumptions for BL models, this measure can not be a product-distribution. (Instead, a shock-like stationary product-measure is described there for certain type of rates, under which the slope of the surface differs on the left side from that on the right side.) 

For SE and some types of ZR processes, law of large numbers \eqref{eq:nsztv} is known. This law and the second moment condition \eqref{eq:masodq} for BL and ZR models possessing convexity condition \ref{con:cvx} are proven in section \ref{sc:nszt}. As shown in the next theorem, this allows us to do further computations on the space-time correlations of the models. We need the following properties of the canonical measure:
\begin{lm}\label{lm:mus}
(i) The sum 
\[
\sum_{z\in I}\sum_{y=z+1}^{\omax}g(y)\frac{\mu(y)}{\sqrt{\mu(z)}}
\]
is convergent, and

\noindent
(ii) the sum
\[
\sum_{z\in I}\sum_{y=z+1}^{\omax}g(y)\,\mu(y)=\Vv(\om_0)
\]
is convergent and the equality holds.
\end{lm}
\begin{proof}
For $\te\in(\un\te,\,\bar\te)$, the tails of the measure $\mu_\te(\cdot)$ have exponential decay. Hence the convergence in both expressions holds. The identity in (ii) is straightforward and is left to the reader.
\end{proof}
The next lemma shows the essential connection of the defect tracer to space-time correlations in the model.
\begin{lm}\label{lm:szorasnegyzet}
Assume condition \ref{con:con} with speed value $C$. Let $B(t)$ be a real-valued function with $\lim_{t\to\infty}B(t)=B\in\Rb$,\ $n_1,\,n_2\in\Zb$,\ $A\in\Rb,\ V_1<V_2$ in $\Rb\cup\{-\infty,\,\infty\}$ and the real interval $\Vf:\,=[V_1,\,V_2]$. If either
\begin{itemize}
\item [(i)] $C\ne V_1,\,V_2$, or
\item [(ii)] $C\in\Rb$ and $A\cdot C=-B$ 
\end{itemize}
holds, then
\begin{multline*}
\lim_{t\to\infty}\sum_{n=\lc tV_1\rc+n_1}^{\lf tV_2\rf+n_2}\left(\frac{n}{t}A+B(t)\right)\cdot\Ev(\widetilde\om_0(0)\,\widetilde\om_n(t))=\\
=(A\,C+B)\cdot{\bf1}\{C\in\Vf\}\cdot\Vv(\om_0),
\end{multline*}
where $\Vv(\om_0)$ is the variance of $\om_0$ w.r.t.\ the canonical Gibbs-measure.
\end{lm}
\begin{proof}
We define $\Vf^t$ by
\[
\Vf^t:\,=\left[V_1+\frac{n_1}{t},\ V_2+\frac{n_2}{t}\right].
\]
By corollary \ref{cor:g-s},
\begin{multline}
\lim_{t\to\infty}\sum_{n=\lc tV_1\rc+n_1}^{\lf tV_2\rf+n_2}\left(\frac{n}{t}A+B(t)\right)\cdot\Ev(\widetilde\om_0(0)\,\widetilde\om_n(t))=\\
=\lim_{t\to\infty}\sum_{n=\lc tV_1\rc+n_1}^{\lf tV_2\rf+n_2}\left(\frac{n}{t}A+B(t)\right)\,\Ev\left({\bf1}\{Q(t)=n\}\cdot\sum_{y=\om_0(0)+1}^{\omax}g(y)\frac{\mu(y)}{\mu(\om_0(0))}\right)=\\
=\lim_{t\to\infty}\Ev\left(\left(A\frac{Q(t)}{t}+B(t)\right)\cdot{\bf1}\{Q(t)/t\in\Vf^t\}\cdot\sum_{y=\om_0(0)+1}^{\omax}g(y)\frac{\mu(y)}{\mu(\om_0(0))}\right)=\\
=\lim_{t\to\infty}\sum_{z\in I}\Ev\left(\left(A\frac{Q(t)}{t}+B(t)\right)\cdot{\bf1}\{Q(t)/t\in\Vf^t\}\cdot{\bf1}\{\om_0(0)=z\}\right)\times\\
\times\sum_{y=z+1}^{\omax}g(y)\frac{\mu(y)}{\mu(z)}.\label{eq:limsum}
\end{multline}
We show that the limit and the summation can be interchanged in this expression. We use Cauchy's inequality to obtain
\begin{multline*}
\left|\Ev\left(\left(A\frac{Q(t)}{t}+B(t)\right)\cdot{\bf1}\{Q(t)/t\in\Vf^t\}\cdot{\bf1}\{\om_0(0)=z\}\right)\right|\le\\
\le\sqrt{\Ev\left(\left(A\frac{Q(t)}{t}+B(t)\right)^2\right)}\cdot\sqrt{\Pv\left(\frac{Q(t)}{t}\in\Vf^t\ \text{and}\ \om_0(0)=z\right)}\le\\
\le K'\cdot\sqrt{\Pv\left(\frac{Q(t)}{t}\in\Vf^t\,\biggl|\,\om_0(0)=z\right)}\cdot\sqrt{\mu(z)}\le K'\cdot\sqrt{\mu(z)}
\end{multline*}
for some constant $K'$ by \eqref{eq:masodq}. Since $g(y)$ is monotone in $y$ and
\[
\sum_{y=\omin}^{\omax}g(y)\,\mu(y)=0,
\]
the sum
\[
\sum_{y=z+1}^{\omax}g(y)\,\mu(y)
\]
is non-negative for any $z\in I$. Hence we can bound from above the absolute value of the terms in \eqref{eq:limsum} for each $z\in I$ by
\[
K'\cdot\sum_{y=z+1}^{\omax}g(y)\frac{\mu(y)}{\sqrt{\mu(z)}},
\]
and the sum 
\[
\sum_{z\in I}K'\cdot\sum_{y=z+1}^{\omax}g(y)\frac{\mu(y)}{\sqrt{\mu(z)}}
\]
is convergent by lemma \ref{lm:mus}. Using dominated convergence, we write
\begin{multline*}
\lim_{t\to\infty}\sum_{n=\lc tV_1\rc+n_1}^{\lf tV_2\rf+n_2}\left(\frac{n}{t}A+B(t)\right)\cdot\Ev(\widetilde\om_0(0)\,\widetilde\om_n(t))=\\
=\sum_{z\in I}\lim_{t\to\infty}\Ev\left(\left(A\frac{Q(t)}{t}+B(t)\right)\cdot{\bf1}\{Q(t)/t\in\Vf^t\}\cdot{\bf1}\{\om_0(0)=z\}\right)\times\\
\times\sum_{y=z+1}^{\omax}g(y)\frac{\mu(y)}{\mu(z)}.
\end{multline*}
We introduce the set $\Vf^t_\ve:\,=\Vf^t\cap\Bf_\ve(C)$, where for $\ve>0,\ \Bf_\ve(C)=(C-\ve,\,C+\ve)\subset\Rb$. Hence $\Vf^t=\Vf^t_\ve\cup\left(\Vf^t\setminus\Bf_\ve(C)\right)$:
\begin{multline}
\lim_{t\to\infty}\sum_{n=\lc tV_1\rc+n_1}^{\lf tV_2\rf+n_2}\left(\frac{n}{t}A+B(t)\right)\cdot\Ev(\widetilde\om_0(0)\,\widetilde\om_n(t))=\\
=\sum_{z\in I}\lim_{t\to\infty}\Ev\left(\left(A\frac{Q(t)}{t}+B(t)\right)\cdot{\bf1}\{Q(t)/t\in\Vf^t_\ve\}\cdot{\bf1}\{\om_0(0)=z\}\right)\times\\
\times\sum_{y=z+1}^{\omax}g(y)\frac{\mu(y)}{\mu(z)}+\\
+\sum_{z\in I}\lim_{t\to\infty}\Ev\left(\left(A\frac{Q(t)}{t}+B(t)\right)\cdot{\bf1}\{Q(t)/t\in\Vf^t\setminus\Bf_\ve(C)\}\cdot{\bf1}\{\om_0(0)=z\}\right)\times\\
\times\sum_{y=z+1}^{\omax}g(y)\frac{\mu(y)}{\mu(z)}.\label{eq:ebont}
\end{multline}
\eqref{eq:ebont} contains two terms. We use Cauchy's inequality on the second term as we have done before:
\begin{multline*}
\left|\Ev\left(\left(A\frac{Q(t)}{t}+B(t)\right)\cdot{\bf1}\{Q(t)/t\in\Vf^t\setminus\Bf_\ve(C)\}\cdot{\bf1}\{\om_0(0)=z\}\right)\right|\le\\
\le K'\cdot\sqrt{\Pv\left(\frac{Q(t)}{t}\in\Vf^t\setminus\Bf_\ve(C)\ \text{and}\ \om_0(0)=z\right)}\le\\
\le K'\cdot\sqrt{\Pv\left(\frac{Q(t)}{t}\notin\Bf_\ve(C)\right)}\to0
\end{multline*}
as $t\to\infty$ by the law of large numbers \eqref{eq:nsztv}. Only the first term of \eqref{eq:ebont} remained, for which we write
\begin{multline}
\lim_{t\to\infty}\sum_{n=\lc tV_1\rc+n_1}^{\lf tV_2\rf+n_2}\left(\frac{n}{t}A+B(t)\right)\cdot\Ev(\widetilde\om_0(0)\,\widetilde\om_n(t))=\\
=\sum_{z\in I}\lim_{t\to\infty}(A\cdot C+B(t)+\mathcal O(\ve))\cdot\Pv\left(Q(t)/t\in\Vf^t_\ve\ \text{and}\ \om_0(0)=z\right)\times\\
\times\sum_{y=z+1}^{\omax}g(y)\frac{\mu(y)}{\mu(z)}.\label{eq:lemveg}
\end{multline}
We have three possibilities.

\noindent
(i) If $C\in\Vf,\ C\ne V_1,\,V_2$, then for small $\ve$ and large $t,\ \Vf^t_\ve=\Bf_\ve(C)$, and by \eqref{eq:nsztv}, 
\begin{multline*}
\lim_{t\to\infty}\Pv\left(Q(t)/t\in\Vf^t_\ve\ \text{and}\ \om_0(0)=z\right)=\\
=\lim_{t\to\infty}\Pv\left(Q(t)/t\in\Bf_\ve(C)\ \text{and}\ \om_0(0)=z\right)=\Pv(\om_0(0)=z)=\mu(z).
\end{multline*}
Hence we can continue \eqref{eq:lemveg} by
\begin{multline*}
\lim_{t\to\infty}\sum_{n=\lc tV_1\rc+n_1}^{\lf tV_2\rf+n_2}\left(\frac{n}{t}A+B(t)\right)\cdot\Ev(\widetilde\om_0(0)\,\widetilde\om_n(t))=\\
=\sum_{z\in I}\lim_{t\to\infty}(A\cdot C+B(t)+\mathcal O(\ve))\cdot\sum_{y=z+1}^{\omax}g(y)\mu(y)\to\\
(A\cdot C+B)\cdot\sum_{z\in I}\sum_{y=z+1}^{\omax}g(y)\mu(y)=(A\cdot C+B)\cdot\Vv(\om_0)
\end{multline*}
as $\ve\to0$. The last equality is a result of lemma \ref{lm:mus}.

\noindent
(ii) If $A\cdot C=-B$, then the right-hand side of \eqref{eq:lemveg} tends to $\mathcal O(\ve)$ as $t\to\infty$ for all $\ve>0$, hence is zero in this limit. Here we used that 
\[
\Pv\left(Q(t)/t\in\Vf^t_\ve\ \text{and}\ \om_0(0)=z\right)\le\Pv(\om_0(0)=z)=\mu(z),
\]
and that
\[
\sum_{z\in I}\sum_{y=z+1}^{\omax}g(y)\mu(y)
\]
is convergent.

\noindent
(iii) In case $C\notin\Vf$, then $\Vf^t_\ve$ is empty for $\ve$ small and $t$ large enough, and hence the right-hand side of \eqref{eq:lemveg} is zero.

\noindent
The result of these three cases completes the proof the lemma.
\end{proof}
Now we are able to compute $\lim_{t\to\infty}\frac{\Vv(J^V(t))}{t}$ for $V=0$. The proof of the general formula \eqref{eq:fovegso} requires some more computations in the next subsection.
\begin{tm}\label{tm:formula2}
Assume condition \ref{con:con} with speed $C$. Then
\begin{eqnarray}
\lim_{t\to\infty}\frac{\Vv(J(t))}{t}&=&\Ev(r)-2\,\Ev(r^*(0)\cdot\widetilde\om_1(0))+2\,C^+\cdot\Vv(\om_0)=\nonumber\\
&=&\Ev(r)+2\,\Ev(r^*(0)\cdot\widetilde\om_0(0))+2\,C^-\cdot\Vv(\om_0).\label{eq:formula2}
\end{eqnarray}
Here $0\le C^\pm$ is the positive or the negative part of $C$, respectively.  
\end{tm}
\begin{proof}
We consider the result of theorem \ref{tm:formula1}. Dividing \eqref{eq:formula1} by $t$ and taking the limit $t\to\infty$ allows us to use the result of lemma \ref{lm:szorasnegyzet}. For the first equality of \eqref{eq:formula1}, we use this lemma with parameters $V_1=0,\ V_2=\infty,\ n_1=1,\ n_2=0,\ A=1,\ B(t)=0$. Then we obtain 
\[
\lim_{t\to\infty}\frac{\Vv J(t)}{t}=\Ev(r)-2\,\Ev(r^*(0)\cdot\widetilde\om_1(0))+2\,C\cdot{\bf 1}\{C\ge0\}\cdot\Vv(\om_0).
\]
For the second equality of \eqref{eq:formula1}, we rewrite the sum as
\[
\lim_{t\to\infty}\frac{\Vv J(t)}{t}=\Ev(r)+2\,\Ev(r^*(0)\cdot\widetilde\om_0(0))-2\sum_{n=-\infty}^{-1}n\,\Ev(\tom_0(0)\,\tom_n(t)),
\]
in order to use lemma \ref{lm:szorasnegyzet} with parameters $V_1=-\infty,\ V_2=0,\ n_1=0,\ n_2=-1,\ A=1,\ B(t)=0$. Hence
\[
\lim_{t\to\infty}\frac{\Vv J(t)}{t}=\Ev(r)+2\,\Ev(r^*(0)\cdot\widetilde\om_0(0))-2\,C\cdot{\bf 1}\{C\le0\}\cdot\Vv(\om_0),
\]
which proves the second equality of the theorem.
\end{proof}

\noindent
We obtained two formulas for the variance of $J(t)$. If the characteristic speed $C$ exists, then we can compute it by subtracting the two lines of (\ref{eq:formula2}).

\fej{The growth in non-vertical directions}\label{sc:sl}

We have examined so far $\Vv(J(t))$, the growth fluctuation of a fixed column, i.e.\ the fluctuation of vertical growth. In this section we deal with the growth fluctuation of the surface in equilibrium, but considered in a slanting direction, namely, $\Vv(J^{(V)}(t))$. From now on, we assume without loss of generity $h_0(0)=0$. 
\begin{proof}[Proof of \eqref{eq:jfnsztv}]
By definition $\om_j(t)=h_{j-1}(t)-h_j(t)$, we have
\begin{equation}
h_i(t)=h_0(t)-\sum_{j=1}^i \om_j(t)\label{eq:h}
\end{equation}
for any site $i>0$, hence for $V>0$,
\[
J^{(V)}(t)=h_\ta(t)=h_0(t)-\sum_{j=1}^\ta \om_j(t)=h_0(t)-\ta\frac{1}{\ta}\sum_{j=1}^\ta \om_j(t).
\]
By ergodicity, the first term has the limit $\Ev r$ a.s.\ when divided by $t$. The second term is $\ta$ times the average of an increasing number of different iid.\ variables. These variables have finite moments, hence the fourth-moment argument (see e.g.\ \cite[Theorem 7.1]{lamp}) is applicable with the discretization series $t_n:\,=n/V$ to show that
\[
\lim_{n\to\infty}\frac{1}{\lf t_nV\rf}\sum_{j=1}^{\lf t_nV\rf}\om_j(t_n)=\lim_{n\to\infty}\frac{1}{n}\sum_{j=1}^n\om_j(t_n)=\Ev(\om)\ \ \text{a.s.}
\]
This shows $\eqref{eq:jfnsztv}$ for the limit taken along the subsequence $t_n$. For any $t\in\Rb^+$, there is a unique $n_t\in\Zb^+$ for which $t_{n_t}\le t<t_{n_t+1}$, and $J^{(V)}(t)-J^{(V)}(t_{n_t})$ is the number of bricks laid on column $n_t$ in a time interval shorter than $1/V$, hence dividing it by $t$ leads a.s.\ to zero in the limit. Therefore $\eqref{eq:jfnsztv}$ holds for the limit of $J^{(V)}(t)/t$ as well. Similar computation works for $V<0$, and finally, the case $V=0$ is trivial.
\end{proof}
Now we consider the fluctuations (with tilde meaning the mean value subtracted).
\begin{multline}
\Vv(J^{(V)}(t))=\Ev\left\{\left(J^{(V)}(t)-\Ev J^{(V)}(t)\right)^2\right\}=\\
=\Ev\left\{\left(\wh_\ta(t)\right)^2\right\}=\Ev\left\{\left(\left[\wh_\ta(t)-\wh_\ta(0)\right]+\wh_\ta(0)\right)^2\right\}=\\
=\Ev\left\{\left(\wh_\ta(t)-\wh_\ta(0)\right)^2\right\}-\Ev\left\{\left(\wh_\ta(0)\right)^2\right\}+2\,\Ev\left(\wh_\ta(t)\,\wh_\ta(0)\right).\label{eq:mozgv}
\end{multline}
By translation-invariance, the first term is $\Vv(J(t))$, computed in the previous sections. By \eqref{eq:h} and by product structure of the measure, the second term of the right-hand side of \eqref{eq:mozgv} is
\[
-\Ev\left\{\left(\wh_\ta(0)\right)^2\right\}=-\ta\cdot\Ev(\tom_0(0)^2)=-\ta\cdot\Vv(\om).
\]
The limit of the third term divided by $t$ in \eqref{eq:mozgv} is computed in the following two lemmas:
\begin{lm}\label{lm:mozkorr}
For $V>0$,
\begin{multline}
\Ev\left(\wh_\ta(t)\,\wh_\ta(0)\right)=\\
=\sum_{n=-\infty}^{\ta-1}(\ta-n)\Ev(\tom_n(t)\,\tom_0(0))+\sum_{n=-\infty}^{-1}n\Ev(\tom_n(t)\,\tom_0(0)).\label{eq:mozkorr}
\end{multline}
\end{lm}
\begin{proof}
Using \eqref{eq:h} again,
\begin{multline}
\Ev\left(\wh_\ta(t)\,\wh_\ta(0)\right)=\\
=-\Ev\left(h_0(t)\sum_{j=1}^\ta\tom_j(0)\right)+\Ev\left(\sum_{i=1}^\ta\tom_i(t)\sum_{j=1}^\ta\tom_j(0)\right)=\\
=-\sum_{j=1}^\ta\Ev(h_0(t)\,\tom_j(0))+\sum_{i=1}^\ta\sum_{j=1}^\ta\Ev(\tom_i(t)\,\tom_j(0)).\label{eq:moz3}
\end{multline}
A martingale
\[
H(t):\,=h_0(t)-\int_0^tr_0(s)\di s
\]
with $H(0)=0$ can be separated in order to show that
\begin{multline*}
\Ev(h_0(t)\,\tom_j(0))=\Ev(H(t)\,\tom_j(0))+\int_0^t\Ev(r_0(s)\,\tom_j(0))\di s=\\
=\int_0^t\Ev(\tr_0(s)\,\tom_j(0))\di s.
\end{multline*}
Now we use an argument very similar to the proof of theorem \ref{tm:formula1}. By lemma \ref{lm:elketto}, the $\text{L}^2$-convergence 
\[
-\lim_{\al\to1}(L\psi_\al)(\un\om)=\tr_0
\]
can be used to replace our integral: for $j\ge1$ we continue by
\begin{multline*}
\Ev(h_0(t)\,\tom_j(0))=\int_0^t\Ev(\tr_0(s)\,\tom_j(0))\di s=-\lim_{\al\to1}\int_0^t\Ev(L\psi_\al(s)\,\tom_j(0))\di s=\\
=-\lim_{\al\to1}\int_0^t\frac{\di}{\di s}\Ev(\psi_\al(s)\,\tom_j(0))\di s=\Ev(\psi_\al(t)\,\tom_j(0))-\Ev(\psi_\al(0)\,\tom_j(0)).
\end{multline*}
Using definition \eqref{eq:fidef} of $\psi_\al$ and product structure of the canonical measure,
\begin{multline*}
\Ev(h_0(t)\,\tom_j(0))=\Ev\left(-\sum_{i=0}^\infty\om_{-i}(t)\,\tom_j(0)\right)-\Ev\left(-\sum_{i=0}^\infty\om_{-i}(0)\,\tom_j(0)\right)=\\
=-\sum_{i=0}^\infty\Ev(\tom_{-i}(t)\,\tom_j(0))=-\sum_{i=-\infty}^0\Ev(\tom_i(t)\,\tom_j(0)).
\end{multline*}
Combining this expression with \eqref{eq:moz3} leads to
\begin{multline*}
\Ev\left(\wh_\ta(t)\,\wh_\ta(0)\right)=\sum_{i=-\infty}^0\sum_{j=1}^\ta\Ev(\tom_i(t)\,\tom_j(0))+\sum_{i=1}^\ta\sum_{j=1}^\ta\Ev(\tom_i(t)\,\tom_j(0))=\\
=\sum_{i=-\infty}^\ta\sum_{j=1}^\ta\Ev(\tom_i(t)\,\tom_j(0))=\sum_{i=-\infty}^\ta\sum_{j=1}^\ta\Ev(\tom_{i-j}(t)\,\tom_0(0))
\end{multline*}
by translation-invariance. Changing the summation indices leads to the proof of the lemma. 
\end{proof}
\begin{lm}\label{lm:mozform2}
Assume condition \ref{con:con}. Then for $V>0$,
\[
\lim_{t\to\infty}\frac{1}{t}\,\Ev\left(\wh_\ta(t)\,\wh_\ta(0)\right)=(V-C^+)^+\cdot\Vv(\om).
\]
\end{lm}
\begin{proof}
We use lemma \ref{lm:szorasnegyzet} for the two terms on the right-hand side of \eqref{eq:mozkorr}. For the first one we set $V_1=-\infty,\ V_2=V,\ n_1=0,\ n_2=-1,\ A=-1,\ B(t)=\ta/t$, while for the second term in \eqref{eq:mozkorr} we put $V_1=-\infty,\ V_2=0,\ n_1=0,\ n_2=-1,\ A=1,\ B(t)=0$. One can easily check that for any $C\in\Rb$ and $V>0$, one of the cases (i) or (ii) of lemma \ref{lm:szorasnegyzet} apply. Consequently, we obtain
\begin{multline*}
\lim_{t\to\infty}\frac{1}{t}\,\Ev\left(\wh_\ta(t)\,\wh_\ta(0)\right)=\\
=\left[(V-C)\cdot{\bf 1}\{C\le V\}+C\cdot{\bf 1}\{C\le0\}\right]\cdot\Vv(\om)=(V-C^+)^+\cdot\Vv(\om).
\end{multline*}
\end{proof}
Now we divide equation \eqref{eq:mozgv} by $t$ and take the limit $t\to\infty$. We use the result of lemma \ref{lm:mozform2} to obtain
\begin{equation}
\lim_{t\to\infty}\frac{\Vv(J^{(V)}(t))}{t}=\lim_{t\to\infty}\frac{\Vv(J(t))}{t}+[2(V-C^+)^+-V]\cdot\Vv(\om)\label{eq:mozpoz}
\end{equation}
for $V>0$. 

For $V<0$, we proceed as we did above with $J^{(V)}$ for positive $V$'s. The only important difference is using $\vp_\al$ instead of $-\psi_\al$ in the proof of lemma \ref{lm:mozkorr}. The result of a similar lemma for $V<0$ is
\begin{multline*}
\Ev(\wh_\tf(t)\,\wh_\tf(0))=\\
=-\sum_{n=\tf+1}^\infty(\tf-n)\Ev(\tom_n(t)\,\tom_0(0))-\sum_{n=1}^\infty n\Ev(\tom_n(t)\,\tom_0(0)).
\end{multline*}
Therefore, lemma \ref{lm:szorasnegyzet} is applicable in a similar way as in lemma \ref{lm:mozform2} above. The result of this application is
\[
\lim_{t\to\infty}\frac{1}{t}\,\Ev\left(\wh_\tf(t)\,\wh_\tf(0)\right)=(V+C^-)^-\cdot\Vv(\om).
\]
Computing $\Vv(J^{(V)})$ for $V<0$ as we did in \eqref{eq:mozgv} leads then to
\begin{equation}
\lim_{t\to\infty}\frac{\Vv(J^{(V)}(t))}{t}=\lim_{t\to\infty}\frac{\Vv(J(t))}{t}+[2(V+C^-)^-+V]\cdot\Vv(\om).\label{eq:mozneg}
\end{equation}

Now, assuming condition \ref{con:con}, we can prove \eqref{eq:fovegso} by the result of theorem \ref{tm:formula2}.
\begin{proof}[Proof of theorem \ref{tm:main}]
All time arguments of our variables for this proof are thought to be zero without mentioning it. By \eqref{eq:rrev},
\begin{multline*}
\Ev\left(r^*\cdot(\tom_0-\tom_1)\right)=\Ev\left(r^*\cdot(\om_0-\om_1)\right)=\\
=\Ev\left(r(\om_0+1,\,\om_1-1)\cdot\frac{\mu(\om_0+1)\,\mu(\om_1-1)}{\mu(\om_0)\,\mu(\om_1)}\cdot(\om_0-\om_1)\right)=\\
\Ev(r\cdot(\om_0-\om_1))-2\Ev(r)=\Ev(r\cdot(\tom_0-\tom_1))-2\Ev(r)=-\Ev(r^*\cdot(\tom_0-\tom_1))-2\Ev(r),
\end{multline*}
we used \eqref{eq:revr} in the last step. Hence we obtain
\[
\Ev\left(r^*\cdot(\tom_0-\tom_1)\right)=-\Ev(r).
\]
We have two formulas for the variance $\Vv(J(t))$ by theorem \ref{tm:formula2}, which are used together with \eqref{eq:mozpoz} and \eqref{eq:mozneg} to obtain
\begin{eqnarray*}
\lim_{t\to\infty}\frac{\Vv(J^{(V)}(t))}{t}&=&\Ev(r)-2\,\Ev(r^*\cdot\widetilde\om_1)+(|V-C|+C)\cdot\Vv(\om)=\nonumber\\
&=&\Ev(r)+2\,\Ev(r^*\cdot\widetilde\om_0)+(|V-C|-C)\cdot\Vv(\om)
.
\end{eqnarray*}
We take the average of these two formulas:
\begin{multline*}
\lim_{t\to\infty}\frac{\Vv(J^{(V)}(t))}{t}=\Ev(r)+\Ev\left(r^*\cdot(\tom_0-\tom_1)\right)+|V-C|\cdot\Vv(\om)=\\
=|V-C|\cdot\Vv(\om).
\end{multline*}
\end{proof}

Now it is easy to prove central limit theorem for $J^{(V)}$.
\begin{proof}[Proof of theorem \ref{tm:cht}]
We introduce the drifted form of $J^{(C)}$ by $i\in\Zb$:
\[
J^{(C)}_i(t):\,=h_{\lf Ct\rf+i}(t)-h_i(0)
\]
for $C\ge0$, and 
\[
J^{(C)}_i(t):\,=h_{\lc Ct\rc+i}(t)-h_i(0)
\]
for $C<0$. Due to translation-invariance, the distribution of this quantity is independent of $i$. Hence by \eqref{eq:fovegso}, for $C\ge0$ and $V\ge0$, the variance of
\[
J^{(C)}_{\ta-\lf Ct\rf}(t)=h_{\ta}(t)-h_{\ta-\lf Ct\rf}(0)=J^{(V)}(t)-h_{\ta-\lf Ct\rf}(0)
\]
is $\mathfrak{o}(t)$ as $t\to\infty$. Thus it follows that we only need central limit theorem for $h_{\ta-\lf Ct\rf}(0)$, which is, by \eqref{eq:h} and by $h_0(0)=0$, the sum of $|\ta-\lf Ct\rf|$ number of iid.\ $\om_i(0)$ variables with finite moments. Hence the theorem follows for $V\ge0,\ C\ge0$. For $V\ge0,\ C<0$,
\[
J^{(C)}_{\ta-\lc Ct\rc}(t)=h_{\ta}(t)-h_{\ta-\lc Ct\rc}(0)=J^{(V)}(t)-h_{\ta-\lc Ct\rc}(0),
\]
here we have (and we only need) central limit theorem for the sum of $|\ta-\lc Ct\rc|$ number of iid.\ $\om_i(0)$ variables, which proves the theorem. Similar argument works for $V<0$ also.
\end{proof}

\fej{The motion of the defect tracer}\label{sc:nszt}

With the help of another type of coupling, with any $n\in\Zb^+$, we prove $L^n$-convergence for $Q(t)/t$ of BL and totally asymmetric ZR models in this section. This coupling only works under convexity condition \ref{con:cvx}, which we assume for the rest of the paper. The idea of the proof is the following: we fix our $(\un\om,\,Q)$ pair and compare it with another model $\un\ze$. The difference between $\un\om$ and $\un\ze$ is realized by second class particles. The current of these particles satisfies law of large numbers by separate ergodicity of $\un\om$ and $\un\ze$, and we compare their motion to our defect tracer $Q$ placed on $\un\om$. The main difficulty is finding the way to couple the defect tracer to the second class particles. As shown later, this coupling can not be made directly; we need to introduce a new process called the $S$-{\sl particles}, a random process defined in terms of the second class particles.

We set $\te_1<\te_2$, then there exists a two dimensional measure $\mu$ on $\Zb\times\Zb$, which has marginals $\mu_{\te_1}$ and $\mu_{\te_2}$, respectively, and for which $\mu(x,\,y)=0$ if $x>y$. We fix two configurations $\un\eta$ and $\un\ze$ of our model, distributed initially according to a product measure with marginals $\Pv(\eta_i(0)=x,\,\ze_i(0)=y)=\mu(x,\,y)$. Therefore, $\un\eta$ is itself in distribution $\un\mu_{\te_1},\ \un\ze$ is in distribution $\un\mu_{\te_2}$, and $\eta_i(0)\le\ze_i(0)$ for each site $i$ is satisfied. According to the basic coupling described in subsection \ref{sc:bas}, $\eta_i(t)\le\ze_i(t)$ holds for all later times $t$, and we have a positive density of second class particles between these two models. The number of these particles at site $i$ is $\ze_i-\eta_i\ge0$. Hence they are initially distributed according to a product measure but, at later times, only the marginal distributions of $\un\eta$ or of $\un\ze$ will possess a product structure. Note that the joint distribution of the processes is translation invariant.

\alfej{The Palm distribution}

For further applications, we want to select ``a typical second class particle''. We do it as follows. We introduce the drifted form of the models: for $k\in\Zb$,
\[
(\tv_k\,\un\eta)_i:\,=\eta_{i+k},\qquad (\tv_k\,\un\ze)_i:\,=\ze_{i+k}.
\]
If $N\in\Zb^+$ is large enough, we choose uniformly one second class particle among the particles present at sites $-N\le i\le N$. We determine the distribution of the values of a function $g$ depending on $(\un\eta,\,\un\ze)$, as seen from the position $k$ of the randomly selected second class particle. For $N$ large enough, the total number 
\[
\sum_{j=-N}^N(\ze_i-\eta_i)
\]
of second class particles at sites $-N\le i\le N$ is positive, and then
\begin{multline*}
\Ev^{(N)}\left(g(\tv_k\,\un\eta,\,\tv_k\,\un\ze)\right)=\Ev\left[\Ev\left(g(\tv_k\,\un\eta,\,\tv_k\,\un\ze)\,\bigr|\,\un\eta,\,\un\ze\right)\right]=\\
=\Ev\left(\sum_{i=-N}^Ng(\tv_i\,\un\eta,\,\tv_i\,\un\ze)\cdot\frac{\ze_i-\eta_i}{\sum_{j=-N}^N(\ze_j-\eta_j)}\right)=\\
=\Ev\left(\frac{\frac{1}{2N+1}\sum_{i=-N}^Ng(\tv_i\,\un\eta,\,\tv_i\,\un\ze)\cdot(\ze_i-\eta_i)}{\frac{1}{2N+1}\sum_{j=-N}^N(\ze_j-\eta_j)}\right).
\end{multline*}
For bounded $g$, the random variable we see in the last line of the display is bounded, and is the quotient of two random variables, both having a.s.\ limit as $N\to\infty$ by translation invariance and ergodicity of translations. Hence our expression converges due to dominated convergence, and have the limit
\begin{equation}
\widehat\Ev\left(g(\un\eta,\,\un\ze)\right):\,=\lim_{N\to\infty}\Ev^{(N)}\left(g(\tv_k\,\un\eta,\,\tv_k\,\un\ze)\right)=\frac{\Ev\left(g(\un\eta,\,\un\ze)\cdot(\ze_0-\eta_0)\right)}{\Ev\left(\ze_0-\eta_0\right)}.\label{eq:palm}
\end{equation}
The distribution $\widehat{\un\mu}$ defined by \eqref{eq:palm} is called the Palm distribution of the process. The Palm measure can be extended to non-negative functions $g$, see \cite{posto}. Note that $\widehat\Pv(\ze_0(0)-\eta_0(0)>0)=1$ according to this measure, i.e.\ we necessarily have at least one second class particle at the origin, if looking the process ``as seen from a typical second class particle''.

By initial product distribution of $(\un\eta,\,\un\ze),\ \widehat{\un\mu}$ is initially also a product measure, consisting of the original marginals $\mu$ for sites $i\ne0$, and of marginal
\begin{equation}
\widehat\mu(x,\,y):\,=\frac{\mu(x,\,y)\cdot(y-x)}{\Ev(\ze_0-\eta_0)}\label{eq:mupal}
\end{equation}
for site $i=0$. For later use, we introduce the pair $(\un\eta'(t),\,\un\ze'(t))$ started from this initial product distribution $\widehat{\un\mu}$.

\alfej{Random walk on the second class particles}

We label the second class particles between $\un\eta$ and $\un\ze$ in space-order. Let $U^{(m)}(t)$ denote the position of the $m$-th second class particle at time $t$. Initially, we look for the first site possessing second class particle on the right side of the origin. We choose one of the particles at this site, giving it label $m=0$:
\[
U^{(0)}(0):\,=\min\left\{i\ge0\,:\,\ze_i>\eta_i\right\}.
\]
We label the particles at $t=0$ in such a way that $U^{(m)}(0)\le U^{(m+1)}(0)\ (\forall m\in\Zb)$ (the order of particles at the same site is not important). We define $J^{(2^{\text{nd}})}_i(t)$ to be the algebraic number of second class particles passing the column between $i$ and $i+1$ in the time interval $[0,\,t]$. This quantity is determined by the evolution of the processes $\un\eta$ and $\un\ze$. For $t=0$, we define 
\begin{equation}
m_i(0):\,=\max\{m\,:\,U^{(m)}(0)\le i\},\label{eq:sorr}
\end{equation}
while for $t>0$,
\[
m_i(t):\,=m_i(0)-J^{(2^{\text{nd}})}_i(t).
\]
We label the particles at later times such that \eqref{eq:sorr} holds at any time $t$ as well. This method assures $U^{(m)}(t)\le U^{(m+1)}(t)$ for all time $t$. The particles labeled from $m_{i-1}+1$ up to $m_i,$ exactly $\ze_i-\eta_i=m_i-m_{i-1}$ of them are at site $i$. (At sites $i$ for which $m_i=m_{i-1}$, there is no second class particle).

We have defined so far the coupled pair $\un\eta$ and $\un\ze$ with the $U^{(m)}(t)$ process of the second class particles indexed in space order at any time $t$. The latter will serve us as a background environment for a new random process, $\left(s^{(n)}(t)\right)_{n\in\Zb}$. Initially, we put $s^{(n)}(0):\,=n$ for each $n$. Assume that just before a second class particle jumps from a site $i$ at a time $t,\ s^{(n)}(t)\in\{m_{i-1}(t)+1,\,m_{i-1}(t)+2,\dots,m_i(t)\}$, which means $U^{(s^{(n)})}(t)=i$ just before the jump. Then by the time $t+0$ of this jump, $s^{(n)}(t+0):\,=\Pi_i(s^{(n)}(t))$, where $\Pi_i$ is a random uniform permutation on the integer set $\{m_{i-1}(t)+1,\,m_{i-1}(t)+2,\dots,m_i(t)\}$.

We can represent this new process as follows. Initially, we put an extra particle, which we call $S$-particle, on each second class particle. The $S$-particles are labeled by $n$, and initially we put the $n$-th $S$-particle on the $n$-th second particle. $s^{(n)}(t)$ stands for the index of the second class particle carrying the $n$-th $S$-particle. Whenever a jump of second class particle happens from site $i$, we permute uniformly and randomly the $S$-particles present at site $i$ just before the jump. According to the labeling of second class particles, one jumping to the right (or to the left, respectively) from site $i$ has index $m_i(t)$ (or $m_{i-1}(t)+1$, respectively) and is carrying exactly the $n$-th $S$-particle, for which $s^{(n)}(t+0)=\Pi_i(s^{(n)}(t))=m_i(t)$ (or $m_{i-1}(t)+1$, respectively). Hence a uniformly and randomly chosen $S$-particle is taken from the site $i$ with the jumping second class particle.

For simplicity, we define $s(t):\,=s^{(0)}(t)$ and $S(t):\,=U^{(s(t))}$, and by simply saying the $S$-particle, we mean the zeroth $S$-particle at site $S(t)$. Then $S(t)$ represents a random walk moving always together with a second class particle, but having always probability $1/(m_i-m_{i-1})=1/(\ze_i-\eta_i)$ of jumping together with a second class particle jumping from the site $i$. As can be derived from table \ref{tab:bas}, the rate for a second class particle to jump to the left (or to the right) is $f(-\eta_i)-f(-\ze_i)$ (or $f(\ze_i)-f(\eta_i)$, respectively). Hence the rate for the $S$-particle to jump to the left (or to the right) together with the jumping second class particle from site $i=S(t)$ is 
\begin{equation}
\frac{f(-\eta_i)-f(-\ze_i)}{\ze_i-\eta_i}\ \  \text{(or}\ \frac{f(\ze_i)-f(\eta_i)}{\ze_i-\eta_i},\ \text{respectively).}\label{eq:srata}
\end{equation}

Recall that $S(0)=U^{(s(0))}(0)=U^{(0)}(0)$ is the first site on the right-hand side of the origin initially with second class particles. We introduce the notation $(\un\eta''(t),\,\un\ze''(t)):\,=(\tv_{S(0)}\un\eta(t),\,\tv_{S(0)}\un\ze(t))$, which is the $(\un\eta(t),\,\un\ze(t))$ process shifted to this initial position $S(0)$ of the $S$-particle. We also introduce its $S''$-particle: $S''(t):\,=S(t)-S(0)$. Hence the initial distribution of $(\un\eta''(0),\,\un\ze''(0))$ is 
%*Ez uj lefele:
modified according to this random shifting-procedure; we show the details in the proof of the next lemma.
%*Idaig uj.

%a product measure, consisting of the original marginals $\mu$ for each site $i\ne0$, and of marginal $\mu_{\text{cond}}$ for the site $i=0$ of the $S''$-particle. The measure $\mu_{\text{cond}}$ is just $\mu$, conditioned on $\{y>x\}$: for $x,\,y\in\Zb$
%\begin{equation}
%\mu_{\text{cond}}(x,\,y)=\Pv(\eta''_0(0)=x,\,\ze''_0(0)=y)={\bf1}\{y>x\}\cdot\frac{\mu(x,\,y)}{\Pv(\ze_0(0)>\eta_0(0))}.\label{eq:mucon}
%\end{equation}

Using the Palm measures, we show that the expected rates for $S$ to jump are bounded in time.
\begin{lm}\label{lm:palm}
Let $n\in\Zb^+,\ k\in\Zb$, and
\begin{equation}
c_i(t):\,=f(\ze_i(t))-f(\eta_i(t))+f(-\eta_i(t))-f(-\ze_i(t))\label{eq:cdef}
\end{equation}
the rate for any second class particle to jump from site $i$. Then
\[
\Ev\left([c_{S(t)}(t)]^n\cdot[\ze_{S(t)}(t)-\eta_{S(t)}(t)]^k\right)\le K(n,\,k)
\]
uniformly in time.
\end{lm}
\begin{proof}
First we consider the pair $(\un\eta'(0),\,\un\ze'(0))$ defined following \eqref{eq:mupal}. As described there, this is in fact the pair $(\un\eta(0),\,\un\ze(0))$ at time $t=0$, as seen from ``a typical second class particle'', or equivalently, as seen from ``a typical $S$-particle''. In this pair, we have at least one second class particle at the origin, which we call $S'$. We let our process $(\un\eta',\,\un\ze')$ evolve, and we follow this ``typical'' $S'$-particle. Started from the Palm-distribution, this tagged $S'$-particle keeps on ``being typical'' (see \cite{posto}), i.e.\ for a function $g$ of the process as seen by $S'$,
\[
\Ev\left(g(\tv_{S'(t)}\un\eta'(t),\,\tv_{S'(t)}\un\ze'(t))\right)=\widehat\Ev\left(g(\un\eta(t),\,\un\ze(t))\right)
\]
with definition \eqref{eq:palm}.

Now we first show the desired result for the $S'$-particle of $(\un\eta',\,\un\ze')$ instead of the $S$-particle of $(\un\eta,\,\un\ze)$. In the previous display, we put the function
\[
g(\un\eta(t),\,\un\ze(t)):\,=[c_0(t)]^n\cdot[\ze_0(t)-\eta_0(t)]^k,
\]
and we denote by $k^+$ the positive part of $k$. We know that $\ze_0(t)-\eta_0(t)\ge1$ holds $\widehat\Pv$-a.s., hence
\begin{multline*}
\Ev\left([c_{S'(t)}(t)]^n\cdot[\ze'_{S'(t)}(t)-\eta'_{S'(t)}(t)]^k\right)=\widehat\Ev\left([c_0(t)]^n\cdot[\ze_0(t)-\eta_0(t)]^k\right)\le\\
\le\widehat\Ev\left([c_0(t)]^n\cdot[\ze_0(t)-\eta_0(t)]^{k^+}\right)=\frac{\Ev\left([c_0(t)]^n\cdot[\ze_0(t)-\eta_0(t)]^{k^++1}\right)}{\Ev(\ze_0(t)-\eta_0(t))}
\end{multline*}
by (37). The function $c_0(t)$ consists of sums of $f(\pm\eta_0(t))$ and $f(\pm\ze_0(t))$, hence the numerator is an $n+k^++1$-order polinom of these functions and of $\ze_0(t)$,\ $\eta_0(t)$. These are all random variables with all moments finite. Therefore, using Cauchy's inequality, the numerator can be bounded from above by products of moments of either $f(\eta_0(t))$ or $f(\ze_0(t))$ or $\eta_0(t)$, or $\ze_0(t)$. The models $\un\eta$ and $\un\ze$ are both separately in their stationary distributions, hence these bounds are constants in time. The denominator is a positive number due to $\te_2>\te_1$ and strict monotonicity of $\Ev_\te(z)$ in $\te$. We see that we found a bound, uniform in time for the function $g$ of $(\un\eta',\,\un\ze')$ as seen from $S'$.

We need to find similar bound for a function $g$ of the original pair $(\un\eta,\,\un\ze)$, as seen from $S$. This is equivalent to finding a bound for $g$ of $(\un\eta'',\,\un\ze'')$ defined above, as seen from $S''$ of this pair. Let us consider 
%*innen uj:
first the initial distribution of $(\un\eta'',\,\un\ze'')$, which we shall call $\un\mu''$. By definition, it is clear that this distribution is the product of the original marginals $\mu$ for sites $i>0$. Fix a $K$ positive integer and two vectors $\un x,\,\un y\in\Zb^\Zb$. For simplicity we introduce the notations
\[
\begin{array}{rcl}
\un\eta_{[a,\,b]}:\,=\left(\eta_a,\,\dots,\,\eta_b\right)&\qquad\text{and}&\qquad\un\ze_{[a,\,b]}:\,=\left(\ze_a,\,\dots,\,\ze_b\right),\\
\un x_{[a,\,b]}:\,=\left(x_a,\,\dots,\,x_b\right)&\qquad\text{and}&\qquad\un y_{[a,\,b]}:\,=\left(y_a,\,\dots,\,y_b\right)
\end{array}
\]
and, where not written, we consider our models at time zero. We break the events according to the initial position $S(0)$ of the $S$-particle in the original pair $(\un\eta,\,\un\ze)$:
\begin{multline*}
\Pv\left(\un\eta''_{[-K,\,0]}=\un x_{[-K,\,0]},\ \un\ze''_{[-K,\,0]}=\un y_{[-K,\,0]}\right)=\\
=\Pv\left(\un\eta_{[S(0)-K,\,S(0)]}=\un x_{[-K,\,0]},\ \un\ze_{[S(0)-K,\,S(0)]}=\un y_{[-K,\,0]}\right)=\\
=\sum_{n=0}^K\Pv\left(\un\eta_{[n-K,\,n]}=\un x_{[-K,\,0]},\ \un\ze_{[n-K,\,n]}=\un y_{[-K,\,0]},\ S(0)=n\right)+\\
+\sum_{n=K+1}^\infty\Pv\left(\un\eta_{[n-K,\,n]}=\un x_{[-K,\,0]},\ \un\ze_{[n-K,\,n]}=\un y_{[-K,\,0]},\ S(0)=n\right)=\\
=\sum_{n=0}^K\Pv\left(\un\eta_{[n-K,\,n]}=\un x_{[-K,\,0]},\ \un\ze_{[n-K,\,n]}=\un y_{[-K,\,0]}\right)\cdot E_n(\un x,\,\un y)+\\
\sum_{n=K+1}^\infty\Pv\left(\un\eta_{[n-K,\,n]}=\un x_{[-K,\,0]},\ \un\ze_{[n-K,\,n]}=\un y_{[-K,\,0]}\right)\cdot E_K(\un x,\,\un y)\cdot\Pv\{F_{n-K}\},
\end{multline*}
where the function $E_n$ of $\un x$ and $\un y$ is an indicator defined by
\[
E_n(\un x,\,\un y):\,={\bf1}\left\{x_{-n}=y_{-n},\,x_{-n+1}=y_{-n+1},\,\dots,\,x_{-1}=y_{-1},\,x_0<y_0\right\},
\]
and the event $F_{n-K}$ is
\[
F_{n-K}:\,=\left\{\eta_0=\ze_0,\,\eta_1=\ze_1,\,\dots,\,\eta_{n-K-1}=\ze_{n-K-1}\right\}.
\]
The last equality follows from the product structure of $\un\mu$ and from the fact that $S(0)$ is the first site to the right of the origin where $\eta_i\ne\ze_i$. Continuing the computation results in
\begin{multline*}
\Pv\left(\un\eta''_{[-K,\,0]}=\un x_{[-K,\,0]},\ \un\ze''_{[-K,\,0]}=\un y_{[-K,\,0]}\right)=\\
=\prod_{i=-K}^0\mu(x_i,\,y_i)\cdot\left[\sum_{n=0}^KE_n(\un x,\,\un y)+E_K(\un x,\,\un y)\cdot\sum_{n=K+1}^\infty\mu\left\{\eta_0=\ze_0\right\}^{n-K}\right]\\
=\prod_{i=-K}^0\mu(x_i,\,y_i)\cdot\left[\sum_{n=0}^KE_n(\un x,\,\un y)+E_K(\un x,\,\un y)\cdot\frac{\mu\left\{\eta_0=\ze_0\right\}}{\mu\left\{\eta_0<\ze_0\right\}}\right]
\end{multline*}
using translation-invariance.

For later purposes, we are interested in the Radon-Nikodym derivative of the distribution $\un\mu''$ of $(\un\eta'',\,\un\ze'')$ w.r.t.\ the Palm distribution $\un{\widehat\mu}$ of $(\un\eta',\,\un\ze')$. Since both have product of marginals $\mu$ for sites $i>0$, we only have to deal with the left part of the origin. Passing to the limit $K\to\infty$, we have
\begin{multline*}
\frac{\di\un\mu''}{\di\un{\widehat\mu}}(\un x,\,\un y)=\lim_{K\to\infty}\frac{\Pv\left(\un\eta''_{[-K,\,0]}=\un x_{[-K,\,0]},\ \un\ze''_{[-K,\,0]}=\un y_{[-K,\,0]}\right)}{\Pv\left(\un\eta'_{[-K,\,0]}=\un x_{[-K,\,0]},\ \un\ze'_{[-K,\,0]}=\un y_{[-K,\,0]}\right)}=\\
=\lim_{K\to\infty}\frac{\prod_{i=-K}^0\mu(x_i,\,y_i)}{\prod_{i=-K}^{-1}\mu(x_i,\,y_i)\widehat\mu(x_0,\,y_0)}\cdot\left[\sum_{n=0}^K E_n(\un x,\,\un y)+E_K(\un x,\,\un y)\cdot\frac{\mu\left\{\eta_0=\ze_0\right\}}{\mu\left\{\eta_0<\ze_0\right\}}\right]=\\
=\frac{\mu(x_0,\,y_0)}{\widehat\mu(x_0,\,y_0)}\cdot\left[\sum_{n=0}^\infty E_n(\un x,\,\un y)+\lim_{K\to\infty}E_K(\un x,\,\un y)\cdot\frac{\mu\left\{\eta_0=\ze_0\right\}}{\mu\left\{\eta_0<\ze_0\right\}}\right]=\\
=\frac{\mu(x_0,\,y_0)}{\widehat\mu(x_0,\,y_0)}\cdot\sum_{n=0}^\infty E_n(\un x,\,\un y)
\end{multline*}
for $\un{\widehat\mu}$-almost all configurations $(\un x,\,\un y)$. Note that the sum on the right-hand side gives exactly the distance between the origin and the first position $i$ to the left of the origin with $x_i\ne y_i$. Hence this sum is finite for $\un{\widehat\mu}$-almost all configurations $(\un x,\,\un y)$.

In view of this result, we can now obtain our estimates. The main idea here is that the pairs $(\un\eta',\,\un\ze')$ and $(\un\eta'',\,\un\ze'')$ only differ in their initial distribution, hence their behavior conditioned on the same initial configuration agree. This is used for obtaining the third expression, and Cauchy's inequality is used for the fourth one below.
\begin{multline*}
\Ev\left([c_{S''(t)}(t)]^n\cdot[\ze_{S''(t)}(t)-\eta_{S''(t)}(t)]^k\right)=\\
=\int_{\widetilde\Omega\cap\{x_0<y_0\}}\Ev\left([c_{S''(t)}(t)]^n\cdot[\ze_{S''(t)}(t)-\eta_{S''(t)}(t)]^k\,|\,\un\eta''(0)=\un x,\,\un\ze''(0)=\un y\right)\times\\
\times\di\un\mu''(\un x,\,\un y)=\\
=\int_{\widetilde\Omega\cap\{x_0<y_0\}}\Ev\left([c_{S'(t)}(t)]^n\cdot[\ze_{S'(t)}(t)-\eta_{S'(t)}(t)]^k\,|\,\un\eta'(0)=\un x,\,\un\ze'(0)=\un y\right)\times\\
\times\frac{\mu(x_0,\,y_0)}{\widehat\mu(x_0,\,y_0)}\cdot\sum_{n=0}^\infty E_n(\un x,\,\un y)\di\un{\widehat\mu}(\un x,\,\un y)\le\\
\le\Biggl[\int_{\widetilde\Omega\cap\{x_0<y_0\}}\left[\Ev\left([c_{S'(t)}(t)]^n\cdot[\ze_{S'(t)}(t)-\eta_{S'(t)}(t)]^k\,|\,\un\eta'(0)=\un x,\,\un\ze'(0)=\un y\right)\right]^2\times\\
\times\di\un{\widehat\mu}(\un x,\,\un y)\Biggr]^\frac12\cdot\Biggl[\int_{\widetilde\Omega\cap\{x_0<y_0\}}\left[\frac{\mu(x_0,\,y_0)}{\widehat\mu(x_0,\,y_0)}\cdot\sum_{n=0}^\infty E_n(\un x,\,\un y)\right]^2\cdot\di\un{\widehat\mu}(\un x,\,\un y)\Biggr]^\frac12\le\\
\le\Biggl[\int_{\widetilde\Omega\cap\{x_0<y_0\}}\Ev\left([c_{S'(t)}(t)]^{2n}\cdot[\ze_{S'(t)}(t)-\eta_{S'(t)}(t)]^{2k}\,|\,\un\eta'(0)=\un x,\,\un\ze'(0)=\un y\right)\times\\
\times\di\un{\widehat\mu}(\un x,\,\un y)\Biggr]^\frac12\cdot\Biggl[\int_{\widetilde\Omega\cap\{x_0<y_0\}}\frac{\mu(x_0,\,y_0)}{\widehat\mu(x_0,\,y_0)}\cdot\left[\sum_{n=0}^\infty E_n(\un x,\,\un y)\right]^2\cdot\di\un\mu(\un x,\,\un y)\Biggr]^\frac12=\\
=\left[\Ev\left([c_{S'(t)}(t)]^{2n}\cdot[\ze_{S'(t)}(t)-\eta_{S'(t)}(t)]^{2k}\right)\right]^\frac12\times\\
\times\left[\int_{\widetilde\Omega\cap\{x_0<y_0\}}\frac{\Ev(\ze_0-\eta_0)}{y_0-x_0}\cdot\left[\sum_{n=0}^\infty E_n(\un x,\,\un y)\right]^2\cdot\di\un\mu(\un x,\,\un y)\right]^\frac12
\end{multline*}
by (38). The first factor of the last display is finite by the first part of the proof. Using the definition of the indicator $E_n$, the second factor can be bounded from above by
\begin{multline*}
\left[\Ev(\ze_0-\eta_0)\right]^\frac12\cdot\left[\int_{\widetilde\Omega\cap\{x_0<y_0\}}\left[\sum_{n=0}^\infty(2n+1)\cdot E_n(\un x,\,\un y)\right]\cdot\di\un\mu(\un x,\,\un y)\right]^\frac12=\\
=\left[\Ev(\ze_0-\eta_0)\right]^\frac12\cdot\left[\sum_{n=0}^\infty(2n+1)\cdot\mu\{\eta_0=\ze_0\}^n\cdot\mu\{\eta_0<\ze_0\}\right]^\frac12
\end{multline*}
using the product property of $\un\mu$, and is again finite.

\end{proof}

Using the rates for the $S$-particle to move, we can prove the following bound for the moments of $S(t)$:
\begin{pr}\label{pr:masods}
For $n\in\Zb^+$,
\begin{equation}
\Ev\left(\frac{|S(t)|^n}{t^n}\right)<K(n)<\infty\label{eq:masods}
\end{equation}
for all large $t$.
\end{pr}
\begin{proof}
For this proof, we denote the jumping rates \eqref{eq:srata} for the $S$-particle by $\rsl$ and $\rsr$, respectively. For $t>0$, we consider the derivative of the quantity above, using these rates:
\begin{multline*}
\frac{\di}{\di t}\Ev\left(\frac{|S(t)|^n}{t^n}\right)=-\frac{n}{t^{n+1}}\Ev\left(|S(t)|^n\right)+\frac{1}{t^n}\frac{\di}{\di t}\Ev\left(|S(t)|^n\right)=-\frac{n}{t^{n+1}}\Ev\left(|S(t)|^n\right)+\\
+\frac{1}{t^n}\Ev\left[\rsl\cdot(|S(t)-1|^n-|S(t)|^n)+\rsr(\cdot|S(t)+1|^n-|S(t)|^n)\right].
\end{multline*}
For $|S(t)|\ge1$, we can bound our expressions:
\[
\frac{\di}{\di t}\Ev\left(\frac{|S(t)|^n}{t^n}\right)\le-\frac{n}{t^{n+1}}\Ev\left(|S(t)|^n\right)+\frac{2^n}{t^n}\Ev\left((\rsr+\rsl)\cdot|S(t)|^{n-1}\right).
\]
We continue by using H\"older's inequality on the right-hand side:
\begin{multline}
\frac{\di}{\di t}\Ev\left(\frac{|S(t)|^n}{t^n}\right)\le\\
\le-\frac{n}{t}\Ev\left(\frac{|S(t)|^n}{t^n}\right)+\frac{2^n}{t}\left\{\Ev\left[(\rsr+\rsl)^n\right]\right\}^{\frac{1}{n}}\cdot\left\{\Ev\left(\frac{|S(t)|^n}{t^n}\right)\right\}^{\frac{n-1}{n}}.\label{eq:derivos}
\end{multline}
Recall that
\[
\rsr(t)+\rsl(t)=c_{S(t)}(t)\cdot[\ze_{S(t)}(t)-\eta_{S(t)}(t)]^{-1},
\]
hence lemma \ref{lm:palm} is applicable with $k=-n$ to show that
\[
\Ev\left[(\rsr(t)+\rsl(t))^n\right]
\]
is bounded in time. Therefore, \eqref{eq:derivos} can be written in the form
\[
\frac{\di}{\di t}\Ev\left(\frac{|S(t)|^n}{t^n}\right)\le-\frac{n}{t}\Ev\left(\frac{|S(t)|^n}{t^n}\right)+\frac{K'(n)}{t}\cdot\left\{\Ev\left(\frac{|S(t)|^n}{t^n}\right)\right\}^{\frac{n-1}{n}}
\]
with some positive constant $K'(n)$. This means that $\Ev(|S(t)|^n/t^n)$ is bounded from above by a solution of the differential equation
\[
\dot y(t)=-\frac{n}{t}\,y(t)+\frac{K'(n)}{t}\cdot y(t)^\frac{n-1}{n}.
\]
Observe that the right-hand side is negative whenever
\[
y(t)>\left(\frac{K'(n)}{n}\right)^n,
\]
hence assuming $y(t_0)<\infty$ for some $t_0>0,\ y(t)$ is bounded (for all $t>t_0$), which gives the proof.
\end{proof}

Now we show law of large numbers for $s(t)$, and then we can show law of large numbers for $S(t)$. For what follows, $\Ev'$ stands for the expectation values according to the distribution of $\un\eta,\ \un\ze,\ \{U^{(m)}\}_{m\in\Zb}$, i.e.\ our background process which determine $m_i(t)$, also. Let $\Ff(t)$ denote the $\si$-field containing all information about these quantities at time $t$. Then $\Ff(t)$ contains all randomness except for the random permutations on $\left(s^{(n)}\right)_{n\in\Zb}$. With \eqref{eq:cdef}, we also introduce the notations
\begin{equation}
\begin{array}{lcl}
C_i(t):&=&(m_i(t)-m_{i-1}(t))^2\cdot c_i(t),\\
p\,(y,\,t):&=&\Pv(s(t)=y\,|\,\Ff(t)),\\
\text{and}&&\\
A_i(t):&=&\max\limits_{m_{i-1}(t)<y\le m_i(t)}p\,(y,\,t)-\min\limits_{m_{i-1}(t)<y\le m_i(t)}p\,(y,\,t)
\end{array}\label{eq:sokdef}
\end{equation}
if $m_i(t)-m_{i-1}(t)>1$, and $A_i(t):\,=0$ otherwise.
\begin{lm}\label{lm:kadas}
\begin{equation}
\frac{\di}{\di t}\Ev(|s(t)|)\le\sqrt{\Ev'\sum_{i=-\infty}^\infty A_i(t)}\cdot\sqrt{\Ev'\sum_{j=-\infty}^\infty A_j(t)\,C_j^2(t)}.\label{eq:kadas}
\end{equation}
\end{lm}
\begin{proof}
We use convention that the empty sum equals zero.
\begin{multline*}
\frac{\di}{\di t}\Ev(|s(t)|)=\lim_{\ve\to0}\frac{\Ev(|s(t+\ve)|)-\Ev(|s(t)|)}{\ve}=\\
=\lim_{\ve\to0}\sum_{z=-\infty}^{\infty}\frac{\Pv(s(t+\ve)=z)\cdot|z|-\Pv(s(t)=z)\cdot|z|}{\ve}=\\
=\lim_{\ve\to0}\Ev'\sum_{i=-\infty}^\infty\sum_{z=m_{i-1}(t)+1}^{m_i(t)}\!\!\!\!\!\frac{\Pv(s(t+\ve)=z\,|\,\Ff(t))\cdot|z|-\Pv(s(t)=z\,|\,\Ff(t))\cdot|z|}{\ve}.
\end{multline*}
We know that uniform random permutation on the indices present at site $S$ happens at each jump of second class particles from $i$ at time $t$. The basic idea is that this permutation makes the probabilities $p\,(y,\,t)=\Pv(s(t)=y\,|\,\Ff(t))$ equalized between $y=m_{i-1}(t)+1\dots m_i(t)$. This jump happens with rate $c_i(t)$ defined in \eqref{eq:cdef}, hence for a site $i$ with at least one second class particle and for $m_{i-1}(t)+1\le z\le m_i(t)$,
\begin{multline*}
\Pv(s(t+\ve)=z\,|\,\Ff(t))=(1-\ve\,c_i(t))\,\Pv(s(t)=z\,|\,\Ff(t))+\\
+\ve\,c_i(t)\sum_{y=m_{i-1}(t)+1}^{m_i(t)}\frac{\Pv(s(t)=y\,|\,\Ff(t))}{m_i(t)-m_{i-1}(t)}+\mathfrak{o}(\ve).
\end{multline*}
Then we obtain
\begin{multline*}
\frac{\di}{\di t}\Ev(|s(t)|)=\Ev'\sum_{i=-\infty}^\infty c_i(t)\sum_{z=m_{i-1}(t)+1}^{m_i(t)}\biggl(\sum_{y=m_{i-1}(t)+1}^{m_i(t)}\frac{p\,(y,\,t)}{m_i(t)-m_{i-1}(t)}\,|z|-\\
-p\,(z,\,t)\cdot|z|\biggr).
\end{multline*}
There exists a $\pi_i$ permutation of the numbers $\{m_{i-1}(t)+1\dots m_i(t)\}$, for which
\[
\sum_{z=m_{i-1}(t)+1}^{m_i(t)}\sum_{y=m_{i-1}(t)+1}^{m_i(t)}\frac{p\,(y,\,t)}{m_i(t)-m_{i-1}(t)}\,|z|\le\sum_{z=m_{i-1}(t)+1}^{m_i(t)}p\,(z,\,t)\cdot|\pi_i(z)|
\]
holds (by permuting higher values of $|z|$ on higher weights), and hence
\begin{multline*}
\frac{\di}{\di t}\Ev(|s(t)|)\le\Ev'\sum_{i=-\infty}^\infty c_i(t)\sum_{z=m_{i-1}(t)+1}^{m_i(t)}p\,(z,\,t)\cdot(|\pi_i(z)|-|z|)=\\
=\Ev'\sum_{i=-\infty}^\infty c_i(t)\sum_{z=m_{i-1}(t)+1}^{m_i(t)}\Bigl(p\,(z,\,t)-\min_{m_{i-1}(t)<y\le m_i(t)}p\,(y,\,t)\Bigr)\cdot(|\pi_i(z)|-|z|)\le\\
\le\Ev'\sum_{i=-\infty}^\infty c_i(t)\sum_{z=m_{i-1}(t)+1}^{m_i(t)}\Bigl(\max_{m_{i-1}(t)<y\le m_i(t)}p\,(y,\,t)-\\
-\min_{m_{i-1}(t)<y\le m_i(t)}p\,(y,\,t)\Bigr)\cdot(m_i(t)-m_{i-1}(t))=\\
=\Ev'\sum_{i\,:\,m_i(t)>m_{i-1}(t)+1}c_i(t)\Bigl(\max_{m_{i-1}(t)<y\le m_i(t)}p\,(y,\,t)-\\
-\min_{m_{i-1}(t)<y\le m_i(t)}p\,(y,\,t)\Bigr)\cdot(m_i(t)-m_{i-1}(t))^2=\Ev'\sum_{i=-\infty}^\infty A_i(t)\,C_i(t)
\end{multline*}
with definitions \eqref{eq:sokdef}. Finally, we use Schwartz and Cauchy's inequality (for simplicity we do not denote time-dependence of the quantities below):
\begin{multline*}
\frac{\di}{\di t}\Ev(|s(t)|)\le\Ev'\sum_{i=-\infty}^\infty A_i\,C_i=\Ev'\sum_{i=-\infty}^\infty \sqrt{A_i}\,\sqrt{A_i}\,C_i\le\\
\le\Ev'\left[\sqrt{\sum_{i=-\infty}^\infty A_i}\cdot\sqrt{\sum_{j=-\infty}^\infty A_j\,C_j^2}\,\right]\le\sqrt{\Ev'\sum_{i=-\infty}^\infty A_i}\cdot\sqrt{\Ev'\sum_{j=-\infty}^\infty A_j\,C_j^2}.
\end{multline*}
\end{proof}
\begin{lm}\label{lm:korl}
The expression
\[
\sqrt{\Ev'\sum_{j=-\infty}^\infty A_j(t)\,C_j^2(t)},
\]
which is the second factor on the right-hand side of \eqref{eq:kadas}, is a bounded function of time.
\end{lm}
\begin{proof}
Due to definitions \eqref{eq:sokdef}, $A_i$ can be bounded from above by
\[
A_i(t)\le\max\limits_{m_{i-1}(t)<y\le m_i(t)}p\,(y,\,t)\le\sum_{y=m_{i-1}(t)+1}^{m_i(t)}p\,(y,\,t)\,=\,\Pv(S(t)=i\,|\,\Ff(t)),
\]
the probability that our $S$-particle is at site $i$. Hence
\begin{multline*}
\sqrt{\Ev'\sum_{j=-\infty}^\infty A_j(t)\,C_j^2(t)}\le\sqrt{\Ev'\sum_{j=-\infty}^\infty\Pv(S(t)=j\,|\,\Ff(t))\cdot C_j^2(t)}=\\
=\sqrt{\Ev'\left[\Ev\left(C_{S(t)}^2(t)\,|\,\Ff(t)\right)\right]}=\sqrt{\Ev\left(C_{S(t)}^2(t)\right)}.
\end{multline*}
The expectation in the last term is bounded in time by lemma \ref{lm:palm} with $n=2,\ k=4$, since
\[
C_i=\left(m_i-m_{i-1}\right)^2\cdot c_i=c_i\cdot\left(\ze_i-\eta_i\right)^2.
\]
\end{proof}
As we know, for any site $i$, the probabilities $p\,(y,\,t)$ can only change by equalizing between $y=m_{i-1}(t)+1\dots m_i(t)$, and the initial distribution is concentrated on $\{s(t=0)=0\}$. Therefore, at every moment $t$, the function $y\to p\,(y,\,t)$ is unimodal. This is clearly the initial situation, and it stays true after each change of this function. By the equalizing property of the $(p\,(y,\,t))_{y\in\Zb}$ process at a jump of second class particle from site $i$,
\[
\max_{m_{i-1}(t)<z\le m_i(t)}p\,(z,\,t)
\]
can never increase. Hence the global maximum $\max_{z\in\Zb}p\,(z,\,t)$ is also a non-increasing function of $t$, and it is bounded as well. Thus its limit exists, which we denote by $P$. It is believed that $P=0$ but we cannot prove this, and this is not necessary for our arguments. 
\begin{lm} Assume $P>0$. Then the set
\[
\left\{x\in\Zb\,:\,p\,(x,\,t)\ge P\right\}
\]
is always contained in the interval $[-1/P,\,1/P]$.
\end{lm}
\begin{proof}
The statement clearly holds initially. For a discrete interval $[x,\,y]$ (with possibly $x=y$ as well), we introduce the block-average
\[
B_{[x,\,y]}(t):\,=\frac{1}{y-x+1}\sum_{z=x}^yp\,(z,\,t),
\]
and we say that $[x,\,y]$ is a {\sl good block}, if $B_{[x,\,y]}(t)\le\min_{z\in[x,\,y]}1/|z|$ (for site $z=0$, we can write 1 instead of $1/|z|$). Any interval is a good block initially. We show this for any time $t$ as well. More precisely, fix $x\le y$, and assume that at a moment $t$, an equalization in the interval $[u,\,v]$ happens:
\[
p(z,\,t+0)=B_{[u,\,v]}(t)
\]
for each $z\in[u,\,v]$. If each finite interval is a good block at $t$, then we show that $[x,\,y]$ is also a good block after this step, at $t+0$. There are four cases.
\begin{itemize}
\item[(i)] If $[u,\,v]$ and $[x,\,y]$ are disjoint or $[u,\,v]\subset[x,\,y]$, then the block-average of $[x,\,y]$ does not change by this step, hence it keeps on being a good block.

\item[(ii)] If $[x,\,y]\subset[u,\,v]$, then $B_{[x,\,y]}(t+0)=B_{[u,\,v]}(t+0)=B_{[u,\,v]}(t)$, and $[u,\,v]$ was a good block at time $t$, hence $[x,\,y]$ is also a good block after this step.

\item[(iii)] In case $[x,\,y]\setminus[u,\,v]\ne\emptyset,\ [u,\,v]\setminus[x,\,y]\ne\emptyset$ and $B_{[u,\,v]}(t)\ge B_{[x,\,y]\setminus[u,\,v]}(t)$ before the step, then 
\[
B_{[u,\,v]}(t+0)=B_{[u,\,v]}(t)\ge B_{[x,\,y]\setminus[u,\,v]}(t)=B_{[x,\,y]\setminus[u,\,v]}(t+0),
\]
hence $B_{[x,\,y]}(t+0)\le B_{[x,\,y]\cup[u,\,v]}(t+0)$. The latter does not change by the step, thus $[x,\,y]\cup[u,\,v]$ keeps on being a good block, which shows that $[x,\,y]$ is also a good block after the step.

\item[(iv)] In case $[x,\,y]\setminus[u,\,v]\ne\emptyset,\ [u,\,v]\setminus[x,\,y]\ne\emptyset$ and $B_{[u,\,v]}(t)<B_{[x,\,y]\setminus[u,\,v]}(t)$ before the step, then by unimodality, $B_{[u,\,v]}(t)\le B_{[u,\,v]\cap[x,\,y]}(t)$, since the function $z\to p\,(z,\,t)$ has no local minimum. This means that $B_{[u,\,v]\cap[x,\,y]}$ does not increase:
\[
B_{[u,\,v]\cap[x,\,y]}(t+0)=B_{[u,\,v]}(t+0)=B_{[u,\,v]}(t)\le B_{[u,\,v]\cap[x,\,y]}(t).
\]
Since $B_{[x,\,y]\setminus[u,\,v]}$ does not change, $B_{[x,\,y]}$ can not increase either, and $[x,\,y]$ was a good block before the step, thus it keeps on being a good block.
\end{itemize}
Applying this result shows the interval containing any single point $z$ to be a good block, i.e.\ $p\,(z,\,t)<1/P$ for $z\notin[-1/P,\,1/P]$, which completes the proof.
\end{proof}
\begin{lm}
Assume $\lim_{t\to\infty}\max_{z\in\Zb}p\,(z,\,t)=P>0$. Then there are $z,\,y$ neighboring sites in the interval $[-1/P-1,\,1/P+1]$ and a time $T>0$, such that the second class particles indexed by $z$ and $y$ cannot be at the same site after $T$: $U^{(z)}(t)\ne U^{(y)}(t)\ (\forall t>T)$. 
\end{lm}
\begin{proof}
Let
\[
A:\,=\left\{z\in\Zb\,:\,\limsup_{t\to\infty}p\,(z,\,t)=P\right\}\ne\emptyset.
\]
By the previous lemma, $A\subset[-1/P,\,1/P]$, and any index $z_{\text{max}}(t)$, for which $p\,(z_{\text{max}}(t),\,t)$ is maximal (and hence larger than or equals to $P$), is also contained in $[-1/P,\,1/P]$ for any $t$. With fixed $P_1<P$ large enough, there exists a moment $T$, such that $p\,(x,\,t)<P_1$ for any $x\notin A$ and for all $t>T$. Hence by $p\,(z_{\text{max}}(t),\,t)\ge P$, all indices $z_{\text{max}}(t)\in A$ for all $t>T$. Let us fix $z\in A$ and $y\notin A$ neighbors, and $y'\notin A$ the other neighbor of $A$. Then infinitely often for $t>T,\ p\,(z,\,t)\ge P>P_1>p\,(y,\,t)$ and $P_1>p\,(y',\,t)$ happens. In this situation, assume that $p\,(z,\,t)$ decreases due to equalization with its neighbors in $A$. Would the result of this step be $p\,(z,\,t)<P$, all indices $z_{\text{max}}(t)$ would be included in this step by unimodality, hence $p\,(z_{\text{max}}(t+0),\,t+0)<P$ would follow, a contradiction. Thus we see that $p\,(z,\,t)\ge P$ can only be violated by an equalization including $y$ or $y'$. If this equalization also includes all indices $z_{\text{max}}(t)$, then the result must be $p\,(y,\,t)\ge P$ or $p\,(y',\,t)\ge P$ by $p\,(z_{\text{max}}(t+0),\,t+0)\ge P$, pulling out at least $P-P_1$ probability from the set $A$. If this step does not include all $z_{\text{max}}$ indices, then it includes indices all with probability at least $P$, hence pulling out at least $(P-P_1)/2$ probability from the set $A$. Since $t>T,\ z_{\text{max}}(t)\in A$, and hence by unimodality, the joint probability of the set $A$ can only decrease. We conclude that assuming equalizing of probabilities between $z\in A$ and $y$ or $y'\notin A$ infinitely often results in decreasing the joint probability of the finite set $A$ infinitely often by a positive constant, which contradicts $P\le p\,(z_{\text{max}}(t),\,t)$ and $z_{\text{max}}(t)\in A$. 
\end{proof}
Now we can prove law of large numbers for the index $s(t)$ of the second class particles carrying our $S$-particle:
\begin{pr}\label{pr:snsztv}
\[
(\forall\delta>0)\ \ \lim_{t\to\infty}\Pv\left(\left|\frac{s(t)}{t}\right|>\delta\right)=0.
\]
\end{pr}
\begin{proof}
By the previous lemma, we see that for $P>0$ there exists a neighboring pair $(z,\,y)\in[-1/P-1,\,1/P+1]$ of second class particles which will never meet after some $T$. After $T$, the process $s(t)$ can not cross such a pair $(z,\,y)$. By translation invariance, it follows a.s.\ that such pairs appear with positive density on $\Zb$ in this case, thus $s(t)$ is bounded a.s.\ and the statement is true. Hence we assume $P=0$ for the rest of the proof. By unimodality,
\begin{multline}
\sum_{i=-\infty}^\infty A_i(t)=\sum_{i\,:\,m_i(t)>m_{i-1}(t)+1}\Bigl(\max_{m_{i-1}(t)<y\le m_i(t)}p\,(y,\,t)-\\
-\min_{m_{i-1}(t)<y\le m_i(t)}p\,(y,\,t)\Bigr)
\le 2\,\max_{z\in\Zb}
p\,(z,\,t).\label{eq:nullas}
\end{multline}
Indirectly let's assume
\[
(\exists\delta>0)\ (\exists K>0)\ (\forall T>0)\ (\exists t>T)\ :\ \Pv\left(\left|\frac{s(t)}{t}\right|>\delta\right)>K.
\]
Then it follows that
\begin{equation}
\Ev(|s(t)|)>K\,\delta\,t\label{eq:ellent}
\end{equation}
for infinitely many and arbitrarily large $t>0$. By \eqref{eq:nullas} and $P=0$,
\[
\sum_{i=-\infty}^\infty A_i(t)\to 0,
\]
thus by dominated convergence theorem
\[
\sqrt{\Ev'\sum_{i=-\infty}^\infty A_i(t)}\to 0.
\]
Hence by lemma \ref{lm:kadas}
\[
\frac{\di}{\di t}\Ev(|s(t)|)\le\sqrt{\Ev'\sum_{i=-\infty}^\infty A_i(t)}\cdot\sqrt{\Ev'\sum_{j=-\infty}^\infty A_j(t)\,C_j^2(t)}\to0
\]
when $t\to\infty$, as 
\[
\sqrt{\Ev'\sum_{j=-\infty}^\infty A_j(t)\,C_j^2(t)}
\]
is bounded by lemma \ref{lm:korl}. That means that 
\[
\frac{\di}{\di t}\Ev(|s(t)|)
\]
tends to zero as $t\to\infty$, which contradicts (\ref{eq:ellent}).
\end{proof}
Now we show the law of large numbers for $S(t)$, the random walk on the background process $\un\eta$ with parameter $\te_1$ and $\un\ze$ with parameter $\te_2$.
\begin{pr}\label{pr:esnsztv}
Let
\begin{equation}
c(\te_1,\,\te_2):\,=2\,\frac{\ch(\te_2)-\ch(\te_1)}{\Ev_{\te_2}(\ze_0)-\Ev_{\te_1}(\eta_0)}\label{eq:c2nd}
\end{equation}
for BL models, and
\begin{equation}
c(\te_1,\,\te_2):\,=\frac{\e{\te_2}-\e{\te_1}}{\Ev_{\te_2}(\ze_0)-\Ev_{\te_1}(\eta_0)}\label{eq:zc2nd}
\end{equation}
for the ZR process. Then for every $\delta>0$
\begin{equation}
\lim_{t\to\infty}\Pv\left(\left|\frac{S(t)}{t}-c(\te_1,\,\te_2)\right|>\delta\right)=0\label{eq:esnsztv}.
\end{equation}
\end{pr}
\begin{proof}
We show the proposition for BL models, the modification for the ZR process is straightforward. By the coupling rules, if a second class particle jumps from $i$ to $i+1$ then the column $g_i$ of $\un\ze$ between sites $i$ and $i+1$ increases by one. If one jumps from $i+1$ to $i$ then the column $h_i$ of $\un\eta$ increases by one. Hence for the current $J^{(2^{\text{nd}})}_i$ of second class particles defined earlier in this subsection,
\[
J^{(2^{\text{nd}})}_i(t)=(g_i(t)-g_i(0))-(h_i(t)-h_i(0)),
\]
i.e.\ it is the difference between the growth of columns $i$ of $\un\ze$ and of $\un\eta$ until time $t$. Due to separate ergodicity of each $\un\ze$ and $\un\eta$, we have law of large numbers for $g_i(t)-g_i(0)$ and for $h_i(t)-h_i(0)$, since each of these models is distributed according to its ergodic stationary measure. Hence with the expectation of the column growth rates, we have
\begin{multline}
\lim_{t\to\infty}\frac{J^{(2^{\text{nd}})}_i(t)}{t}=\Ev_{\te_2}\left(f(\ze_0)+f(-\ze_0)\right)-\Ev_{\te_1}\left(f(\eta_0)+f(-\eta_0)\right)=\\
=2\,\left(\ch(\te_2)-\ch(\te_1)\right)\ \ a.s.\label{eq:jnsztv}
\end{multline}
We extend definition \eqref{eq:sorr} for $x\in\Rb$:
\[
m_x(t):\,=m_{\lf x\rf}(t)=\max\{m\ :\ U^{(m)}(t)\le x\}
\]
Obviously, $m_x(t)=m_x(0)-J_{\lfloor x\rfloor}^{2^{\text{nd}}}(t)$. If $K\in\Rb$ then
\[
\lim_{v\to\infty}\frac{m_{(Kv)}(0)}{K\,v}=\Ev(\ze_0)-\Ev(\eta_0)=\,:p\ \ \text{a.s.}
\]
since at $t=0$, the starting distribution of the number of second class particles at different sites is a product measure.
\begin{multline}
\lim_{t\to\infty}\Pv\left(\left|\frac{m_{(Kt)}(t)}{t}+(c(\te_1,\,\te_2)-K)p\right|>\ve\right)=\\
=\lim_{t\to\infty}\Pv\left(\left|\frac{m_{(Kt)}(t)}{t}-Kp+2\ch(\te_2)-2\ch(\te_1)\right|>\ve\right)=\\
=\lim_{t\to\infty}\Pv\left(\left|\frac{m_{(Kt)}(0)-J^{(2^{\text{nd}})}_{\lfloor Kt\rfloor}(t)}{t}-Kp+2\ch(\te_2)-2\ch(\te_1)\right|>\ve\right)=\\
=\lim_{t\to\infty}\Pv\left(\left|-\frac{J^{(2^{\text{nd}})}_{\lfloor Kt\rfloor}(t)}{t}+2\ch(\te_2)-2\ch(\te_1)\right|>\ve\right)=\\
=\lim_{t\to\infty}\Pv\left(\left|-\frac{J^{(2^{\text{nd}})}_0(t)}{t}+2\ch(\te_2)-2\ch(\te_1)\right|>\ve\right)=0\label{eq:segedtv}
\end{multline}
by translation-invariance and by \eqref{eq:jnsztv}, for any $\ve>0$. Recall that $S(t)$ is the position of the zeroth $S$-particle, i.e.\ the position of the $s(t)$-th second class particle: $S(t)=U^{(s(t))}(t)$. Hence
\begin{multline}
\Pv\left(\left|\frac{S(t)}{t}-c(\te_1,\,\te_2)\right|>\delta\right)=\Pv\left(\left|\frac{U^{(s(t))}(t)}{t}-c(\te_1,\,\te_2)\right|>\delta\right)=\\
=\Pv\left(U^{(s(t))}(t)>c(\te_1,\,\te_2)\,t+\delta\,t\right)+\Pv\left(U^{(s(t))}(t)<c(\te_1,\,\te_2)\,t-\delta\,t\right).\label{eq:soso}
\end{multline}
In case
\[
U^{(s(t))}(t)>c(\te_1,\,\te_2)\,t+\delta\,t
\]
it follows by definitions of $m_x(t)$ and of $p$ that
\begin{multline*}
s(t)>m_{(c(\te_1,\,\te_2)\,t+\delta\,t)}(t)\ \ \text{, hence}\\
\Pv\left(U^{(s(t))}(t)>c(\te_1,\,\te_2)\,t+\delta\,t\right)\le\Pv\left(\frac{s(t)}{t}>\frac{m_{(c(\te_1,\,\te_2)\,t+\delta\,t)}(t)}{t}\right)\le\\
\le\Pv\left(\frac{s(t)}{t}>\frac{\de}{2}\,p\right)+\Pv\left(\frac{m_{(c(\te_1,\,\te_2)\,t+\delta\,t)}(t)}{t}<\frac{\de}{2}\,p\right).
\end{multline*}
As time goes on, the first term goes to zero due to proposition \ref{pr:snsztv}, and so does the second term by \eqref{eq:segedtv} (with $K=c(\te_1,\,\te_2)+\de$).

\noindent
In case
\[
U^{(s(t))}(t)<c(\te_1,\,\te_2)\,t-\delta\,t
\]
it follows that
\begin{multline*}
s(t)\le m_{(c(\te_1,\,\te_2)\,t-\delta\,t)}(t)\ \ \text{, hence}\\
\Pv\left(U^{(s(t))}(t)<c(\te_1,\,\te_2)\,t-\delta\,t\right)\le\Pv\left(\frac{s(t)}{t}\le\frac{m_{(c(\te_1,\,\te_2)\,t-\delta\,t)}(t)}{t}\right)\le\\
\Pv\left(\frac{s(t)}{t}\le-\frac{\de}{2}\,p\right)+\Pv\left(\frac{m_{(c(\te_1,\,\te_2)\,t-\delta\,t)}(t)}{t}>-\frac{\de}{2}\,p\right).
\end{multline*}
The first term again goes to zero due to proposition \ref{pr:snsztv}, and so does the second term by \eqref{eq:segedtv} (with $K=c(\te_1,\,\te_2)-\de)$. Thus we see that both terms on the right-hand side of \eqref{eq:soso} tend to zero as $t\to\infty$.
\end{proof}

\alfej{Coupling the defect tracer to the $S$-particles}

We fix the model $\un\om$ in stationary distribution $\un\mu_\te$ with the defect tracer $Q(t)$ started from the origin. We prove theorem \ref{tm:con} for BL and (totally asymmetric) ZR models. A natural idea would be to couple the defect tracer $Q$ to the second class particles, present at the same site $Q$. The problem is that, either to the left or to the right, the rate for any jump of second class particles form the site $Q$ may be higher than the rate for $Q$ to jump. On the other hand, one second class particle always stays at site $i$ after one jump from $i$, in case more than one of them were present at $i$. The solution is to couple the defect tracer to the $S$-particle, for which the desired conditions are already proven by propositions \ref{pr:masods} and \ref{pr:esnsztv}. For simplicity reasons, in case of the ZR process we let $f(-z):\,=0$ for $z>0$, and hence $\mu(-z)$ of ZR is also zero in these cases.

\bigskip\noindent
{\bf The upper bound for $Q$.}

\bigskip
First, we identify $\un\eta$ distributed according to $\un\mu_{\te_1}$ with $\un\om$ possessing the defect tracer $Q$, therefore we set $\te_1:\,=\te<\te_2$. We have then $\om_i(t)\le\ze_i(t)$ for all $t$ according to the basic coupling, and recall that $Q(0)=0\le S(0)$. In what follows, we are going to couple the random permutations of the $S$-particles, thus the random walk $S(t)$ of the zeroth $S$-particle, with the defect tracer $Q(t)$. We only couple them in case $Q(t)=S(t)$. The basic observation we use is that the rates \eqref{eq:srata} for the jump of the $S$-particle can be compared to the rates for the jump of the defect tracer $Q(t)$. As we have seen at the introduction of BL models, it is enough to consider the ``effect of bricklayers'' standing at each position $i$. That is to say, we are allowed to consider the $\om_i$-dependent parts of $r(\om_{i-1},\,\om_i)$ and $r(\om_i,\,\om_{i+1})$ only, since the $\om_i$-dependent parts are added to the $\om_{i-1}$-dependent or to the $\om_{i+1}$-dependent parts in these rates. In the rest of the paper, we describe couplings by giving rates of bricklayers standing at each site $i$. This observation also holds for the zero range process (by saying rate for a particle to jump instead of saying rate for bricklayers to lay bricks).

In tables \ref{tab:SQl} and \ref{tab:SQr}, $h_i\uparrow$ means that the column of the model $\un\om$ between $i$ and $i+1$ has increased by one, $g_i\uparrow$ means that this column of $\un\ze$ has increased by one, $\curvearrowright$ means the jump to the right from $i$, $\curvearrowleft$ means the jump to the left from $i$.

Note that by $i=S=Q,\ \ze_i\ge\om_i+1$. The rates are non negative due to monotonicity of $f$ and convexity condition \ref{con:cvx}. By summing the rates corresponding to any column of the tables, one can verify that each $\un\om$ and $\un\ze$ evolves according its original rates, $Q$ has the jump rates according to the basic coupling described in table \ref{tab:bas}, and $S$ also has the appropriate rates \eqref{eq:srata}. We see that once being at position $S$, the defect tracer $Q$ can't move right without moving $S$ with it and $S$ can't move left without moving $Q$ with it. Hence our rules preserve the condition $Q\le S$.

We have so far the upper bound $Q(t)\le S(t)$, and we have the law of large numbers \eqref{eq:esnsztv} with speed $c(\te,\,\te_2)$ defined in either \eqref{eq:c2nd} or in \eqref{eq:zc2nd} for any $\te_2>\te$, and the $n$-th moment condition \eqref{eq:masods} for this $S(t)$ process.

\bigskip\noindent
{\bf The lower bound for $Q$.}

\bigskip
Now we show a similar coupling which results in a lower bound for $Q$. The natural idea would be to identify $\un\ze$ with $\un\om$, and couple $Q$ to the $S$-particle. The rates for $Q$ and $S$ to jump with would allow $Q(t)\ge S(t)$. However, this coupling can not be realized in a similar way that the coupling described above: there is no way for $Q$ and $S$ to step together, since only one brick can be laid at a time to a column.

Therefore, we need to modify the initial distribution of the models as follows. Let $\mu(x,\,y)$ be, as before, a two dimensional distribution giving probability zero to $x>y$, and having marginals $\mu_{\te_1}$ and $\mu_{\te_2}$, respectively. Fix the pair $(\un\eta,\,\un\ze)$, as before, with the product of $\mu(x,\,y)$ for different sites as initial distribution. Define
\begin{equation}
\mu'(y,\,x):\,=\mu(x,\,y)\cdot\frac{\mu_{\te_2}(y-1)}{\mu_{\te_2}(y)}.\label{eq:muve}
\end{equation}
Fix the pair $(\un\eta',\,\un\ze')$, with the product of $\mu(x,\,y)$ for each site $i\ne0$ and of $\mu'(x,\,y)$ for the site $i=0$ as initial distribution. Then $\eta'_i(0)\le\ze'_i(0)$ holds a.s.\ for each site $i$, hence the basic coupling is applicable for this pair of models. We have second class particles between $\un\eta'$ and $\un\ze'$, and we introduce the $S'$-particles as well, starting $S'_0$ from the first site on the left-hand side of the origin:
\[
S'(0)=S'_0(0):\,=\max\{i\le0\,:\,\ze'_i(0)>\eta'_i(0)\}.
\]
Assume now that the $S$-particle of $(\un\eta,\,\un\ze)$ is also started from the first site on the left-hand side of the origin, instead of starting it from the right-hand side of the origin:
\[
S(0)=S_0(0):\,=\max\{i\le0\,:\,\ze_i(0)>\eta_i(0)\}.
\]
Then it is clear, that propositions \ref{pr:masods} and \ref{pr:esnsztv} also hold for this $S$-particle. Now we derive these statements for $S'$ as well. Since initially $(\un\eta',\,\un\ze')$ only differs from $(\un\eta,\,\un\ze)$ by the distribution at the origin, the conditional expectations
\begin{equation}
\Ev\left(S'(t)\,|\,\eta'_0(0)=x,\ \ze'_0(0)=y)\right)=\Ev\left(S(t)\,|\,\eta_0(0)=x,\ \ze_0(0)=y)\right)\label{eq:sfelt}
\end{equation}
agree. This is the basic idea of the following
\begin{lm}
The moment condition \eqref{eq:masods} and the law of large numbers \eqref{eq:esnsztv} hold for $S'$ as well.
\end{lm}
\begin{proof}
By the use of \eqref{eq:sfelt} and Cauchy's inequality in a similar way than in the proof of lemma \ref{lm:palm},
\begin{multline*}
\Ev\left(\frac{|S'(t)|^n}{t^n}\right)=\sum_{x\le y}\Ev\left(\frac{|S'(t)|^n}{t^n}\,\Bigr|\,\eta'_0(0)=x,\,\ze'_0(0)=y\right)\cdot\mu'(x,\,y)=\\
=\sum_{x\le y}\Ev\left(\frac{|S(t)|^n}{t^n}\,\Bigr|\,\eta_0(0)=x,\,\ze_0(0)=y\right)\cdot\sqrt{\mu(x,\,y)}\cdot\frac{\mu'(x,\,y)}{\sqrt{\mu(x,\,y)}}\le\\
\le\left[\sum_{x\le y}\left(\Ev\left(\frac{|S(t)|^n}{t^n}\,\Bigr|\,\eta_0(0)=x,\,\ze_0(0)=y\right)\right)^2\cdot\mu(x,\,y)\right]^\frac12\times\\
\times\left[\sum_{x\le y}\frac{\mu'(x,\,y)}{\mu(x,\,y)}\cdot\mu'(x,\,y)\right]^\frac12\le\left[\Ev\left(\frac{S(t)^{2n}}{t^{2n}}\right)\right]^\frac12\cdot\left[\sum_{x\le y}\frac{\mu'(x,\,y)}{\mu(x,\,y)}\cdot\mu'(x,\,y)\right]^\frac12.
\end{multline*}
The first factor of the display is bounded by proposition \ref{pr:masods}. For the second factor, by \eqref{eq:muve} and \eqref{eq:om} we write
\begin{multline*}
\sum_{x\le y}\frac{\mu'(x,\,y)}{\mu(x,\,y)}\cdot\mu'(x,\,y)=\sum_{x\le y}\frac{\mu_{\te_2}(y-1)}{\mu_{\te_2}(y)}\cdot\mu'(x,\,y)=\\
=\sum_{y\in\Zb}\frac{\mu_{\te_2}(y-1)}{\mu_{\te_2}(y)}\cdot\mu_{\te_2}(y-1)=\sum_{y\in\Zb}\frac{f(y)}{\e{\te_2}}\cdot\frac{\e{\te_2(y-1)}}{f(y-1)!}\cdot\frac{1}{Z(\te_2)}=\frac{1}{\e{2\te_2}}\Ev_{\te_2}(f(y)^2),
\end{multline*}
which is again finite. Hence \eqref{eq:masods} holds for $S'$ as well.

For the law of large numbers, we know that for any $\de>0$,
\begin{multline*}
0=\lim_{t\to\infty}\Pv\left(\left|\frac{S(t)}{t}-c(\te_1,\,\te_2)\right|>\delta\right)=\\
=\lim_{t\to\infty}\sum_{x\le y}\Pv\left(\left|\frac{S(t)}{t}-c(\te_1,\,\te_2)\right|>\delta\,\Bigr|\,\eta_0(0)=x,\,\ze_0(0)=y\right)\cdot\mu(x,\,y)=\\
=\lim_{t\to\infty}\sum_{x\le y}\Pv\left(\left|\frac{S'(t)}{t}-c(\te_1,\,\te_2)\right|>\delta\,\Bigr|\,\eta'_0(0)=x,\,\ze'_0(0)=y\right)\cdot\mu(x,\,y),
\end{multline*}
hence \eqref{eq:esnsztv} follows for $S'$ as well by absolute continuity of $\mu'$ w.r.t.\ $\mu$.
\end{proof}

In order to obtain lower bound for $Q$ of $\un\om$ distributed according to $\un\mu_\te$, set $\te_2=\te>\te_1$. The marginal distribution of $\ze_0(0)$ is the second marginal of $\mu'$, namely, $\mu_{\te_2}(y-1)=\mu_\te(y-1)$. Hence it is possible to fix the pair $(\un\eta',\,\un\ze')$ defined above with
\[
\un\ze'(t)=\un\om(t)+\un\de_{Q(t)},\qquad Q(0)=0,
\]
i.e.\ $\un\om$ is coupled to $\un\ze'$ with the defect tracer $Q$ between them. Note that $S'(0)\le0=Q(0)$. We show the coupling that preserves $S'(t)\le Q(t)$ for all later times. We only couple $Q$ to the random permutations acting on $S'$ in case $Q=S'$ for a site $i$. For tables \ref{tab:SQ'l} and \ref{tab:SQ'r}, $h'_i\uparrow$ means that the column of the model $\un\eta'$ between $i$ and $i+1$ has increased by one, $g'_i\uparrow$ means that this column of $\un\ze'$ has increased by one. Note that by $i=S'=Q,\ \ze'_i\ge\eta'_i+1$. As at the coupling for the upper bound, the rates are non negative due to monotonicity of $f$ and convexity condition \ref{con:cvx}. By summing the rates corresponding to any column of the tables, one can verify that each $\un\eta'$ and $\un\ze'$ evolves according its original rates, $Q$ has the jump rates according to the basic coupling described in table \ref{tab:bas} (hence $\un\om$ also evolves according its original rates), and $S'$ also has the appropriate rates \eqref{eq:srata}. We see that once being at position $S'$, the defect tracer $Q$ can't move left without moving $S'$ with it and $S'$ can't move right without moving $Q$ with it. Hence our rules preserve the condition $Q\ge S'$.

\begin{proof}[Proof of theorem \ref{tm:con}] By the upper bound and the lower bound above, we have 
\[
S(t)\ge Q(t)\ge S'(t)
\]
and for any $\te_2>\te>\te_1$, we have weak law of large numbers for $S$ with $c(\te,\,\te_2)$, and for $S'$ with $c(\te_1,\,\te)$, respectively. Hence taking the limits $\te_1\nearrow\te$ and $\te_2\searrow\te$ completes the proof of the law of large numbers \eqref{eq:nsztv} by computing
\[
\lim_{\te_1\nearrow\te}c(\te_1,\,\te)=\lim_{\te_2\searrow\te}c(\te,\,\te_2)=C(\te)
\]
both for BL and ZR models. Moreover, for any $n\in\Zb^+$, we have $n$-th moment condition \eqref{eq:masods} for both $S$ and $S'$, hence not only \eqref{eq:masodq}, but the $n$-th moment condition follows as well for $Q$. This also shows $L^n$-convergence of $Q(t)/t$ for any $n\in\Zb^+$.
\end{proof}

\alfej{Strict monotonicity of $C(\te)$}\label{sc:strict}

As a consequence of the type of coupling methods shown above, we are able to show strict convexity of the function $\mathcal H(\vr)$ of \eqref{eq:haro}. First we refer to the coupling which shows (non strict) convexity, and then we complete the proof of strict convexity by some analytic arguments.
\begin{rem}\label{rm:dupla}
Let $\un\om,\,\un\om'$ be two copies of a model (either BL or ZR model) possessing condition \ref{con:cvx}, with the defect tracers $Q(t)$ and $Q'(t)$, respectively. Assume that for each site $i$ and for time $t=0$ 
\[
\om_i(0)\le\om'_i(0)\ \ \text{and}\ \ Q(0)\le Q'(0).
\]
Then it is possible to couple such way that for all $t\ge0$ and any $i\in\Zb$,
\[
\om_i(t)\le\om'_i(t)\ \ \text{and}\ \ Q(t)\le Q'(t)\ \ \text{a.s.}
\]
is satisfied.
\end{rem}
This coupling is very similar to the ones shown in this subsection, we do not give the details here. The pair $(\un\om,\,\un\om')$ is coupled according to the basic coupling, and we can apply this proposition for the case when their joint distribution has marginals $\un\mu_\te$ and $\un\mu_{\te'}$, respectively. Then we simply see that the motion of the defect tracer of a model has a monotonicity in the parameter $\te$ of the model's stationary distribution. In the introduction we saw that this implies convexity of the function $\mathcal H(\vr)$. We prove now strict convexity of this function:
\begin{proof}[Proof of proposition \ref{pr:cvx}]
First note that by the form \eqref{eq:om} of the measure $\mu_\te$, we have
\[
\begin{array}{rcl}
\vr(\te)=&\Ev_\te(\om)&=\frac{\di}{\di\te}\log\left(Z(\te)\right),\\
\Ev_\te\left({\widetilde\om}^2\right)=&\frac{\di}{\di\te}\Ev_\te\left(\om\right)&>0,\\
\Ev_\te\left({\widetilde\om}^3\right)=&\frac{\di}{\di\te}\Ev_\te\left({\widetilde\om}^2\right)&=\frac{\di}{\di\te}\left(\Vv_\te(\om)\right),
\end{array}
\]
where tilde stands for the centered variable. For the BL model, we need to show strict convexity of the function
\[
\mathcal H(\vr)=\Ev_{\te(\vr)}(r)=\e{\te(\vr)}+\e{-\te(\vr)}.
\]
We compute its derivative
\[
\frac{\di}{\di\vr}\mathcal H(\vr)=\frac{\frac{\di}{\di\te}\left(\e{\te}+\e{-\te}\right)}{\frac{\di\vr}{\di\te}}=\frac{\left(\e{\te}-\e{-\te}\right)}{\Ev_\te\left({\widetilde\om}^2\right)},
\]
and, similarly, the second derivative
\[
\frac{\di^2}{\di\vr^2}\mathcal H(\vr)=\frac{1}{\left[\Ev_\te\left({\widetilde\om}^2\right)\right]^3}\left[\left(\e{\te}+\e{-\te}\right)\,\Ev_\te\left({\widetilde\om}^2\right)-\left(\e{\te}-\e{-\te}\right)\,\Ev_\te\left({\widetilde\om}^3\right)\right].
\]
Hence (strict) positivity of 
\begin{equation}
\left[\left(\e{\te}+\e{-\te}\right)\,\Ev_\te\left({\widetilde\om}^2\right)-\left(\e{\te}-\e{-\te}\right)\,\Ev_\te\left({\widetilde\om}^3\right)\right]\label{eq:anal}
\end{equation}
on an interval of $\te$ is equivalent to (strict) convexity of $\mathcal H(\vr)$ on the corresponding interval of $\vr(\te)$. \eqref{eq:anal} contains derivatives of $\log\left(Z\left(\te\right)\right)$, which is by definition analytic, hence \eqref{eq:anal} is also an analytic function of $\te$. Moreover, by the previous remark, we know convexity of $\mathcal H(\vr)$, hence non-negativity of \eqref{eq:anal}. Since this function is strictly positive at $\te=0$ by symmetry properties of $\mu_\te$, there are at most countably many isolated points at which this analytic function is not strictly positive, hence we have at most countably many isolated points at which the second derivative of $\mathcal H(\vr)$ is not strictly positive. This completes the proof for the BL models.

As for the ZR process, similar computation leads to
\[
\left[\e{\te}\,\Ev_\te\left({\widetilde\om}^2\right)-\e{\te}\,\Ev_\te\left({\widetilde\om}^3\right)\right]
\]
in place of \eqref{eq:anal}. As we know non-negativity of this function by convexity of $\mathcal H(\vr)$, we only need to show $\Ev_\te\left({\widetilde\om}^2\right)\ne\Ev_\te\left({\widetilde\om}^3\right)$ for some $\te$, then the previous analytic argument leads to strict convexity.

Indirectly, assume 
\begin{equation}
\Ev_\te\left({\widetilde\om}^2\right)=\Ev_\te\left({\widetilde\om}^3\right)\label{eq:23}
\end{equation}
for all $\te<\bar\te$. Since the right-hand side is the derivative of the left-hand side, it follows that
\[
\Ev_\te\left({\widetilde\om}^2\right)=A\cdot\e{\te}
\]
for some $A>0$. Integrating this we have
\[
\Ev_\te(\om)=A\cdot\e{\te}
\]
(the additive constant is zero as can be seen by taking the limit $\te\to-\infty$). Integrating again we have
\begin{eqnarray*}
\log\left(Z(\te)\right)=&A\cdot\e{\te}+K,\qquad\text{i.e.}\\
Z(\te)=&K'\cdot\e{A\cdot\e{\te}},\qquad\text{i.e.}\\
\sum_{z=0}^\infty\frac{\e{\te z}}{f(z)!}=&\displaystyle{K'\cdot\sum_{z=0}^\infty\frac{A^z\cdot\e{\te z}}{z!}}
\end{eqnarray*}
for all $\te<\bar\te$, which leads to $f(z)!=z!/A^z,\ f(z)=z/A$. Hence we see that if at least for one $z\ge1$ value we have $f(z+1)-f(z)>f(z)-f(z-1)$, then \eqref{eq:23} is not true for some $\te$, and then strict convexity of $\mathcal H(\vr)$ holds. We also see linearity of $\mathcal H(\vr)$ when $f$ is linear.
\end{proof}

\renewcommand\thesection{}\section{\!\!\!\!\!\!Acknowledgement}

\noindent
The author is grateful to B\'alint T\'oth for initiating the study of bricklayers' models, for giving some of the basic ideas of this paper, and for helping him in many questions. Correcting some mistakes is also acknowledged to him, as well as to Benedek Valk\'o. The author is also grateful to Pablo A.\ Ferrari and to Frank H.\ J.\ Redig for extremely helpful comments on some of the proofs. 

This work was partially supported by the Hungarian National Scientific Research Fund, grant no.\ OTKA T037685 and by the Hungarian Academy of Sciences, TKI Stochastics@TUB.

\bibliography{eredeti}
\bibliographystyle{elsart-num}
\bigskip\bigskip\bigskip\bigskip
\begin{flushright}
{\renewcommand\arraystretch{0.8}
{\small\begin{tabular}{l}
M\'arton Bal\'azs\\
Institute of Mathematics,\\
Technical University Budapest\\
1111. Egry J\'ozsef u.\ 1. H \'ep.\ V. 7.\\
Budapest, Hungary\\
\texttt{balazs@math.bme.hu}
\end{tabular}}}
\end{flushright}

\begin{figure}[p]
\begin{center}
\begin{picture}(100, 120)(0, -10)
\linethickness{0.2pt}
\put(0, 10){\line(1, 0){100}}
\put(0, 30){\line(1, 0){100}}
\put(0, 50){\line(1, 0){100}}
\put(0, 70){\line(1, 0){30}}
\put(90, 70){\line(1, 0){10}}
\put(0, 90){\line(1, 0){10}}

\put(10, 0){\line(0, 1){90}}
\put(30, 0){\line(0, 1){70}}
\put(50, 0){\line(0, 1){50}}
\put(70, 0){\line(0, 1){50}}
\put(90, 0){\line(0, 1){50}}

\put(29, -7){$\scriptstyle{i}$}
\put(45, -7){$\scriptstyle{i\!+\!1}$}

\linethickness{2pt}

\put(0, 110){\line(1, 0){10}}
\put(10, 90){\line(1, 0){20}}
\put(30, 70){\line(1, 0){20}}
\put(50, 50){\line(1, 0){40}}
\put(90, 90){\line(1, 0){10}}

\put(10, 90){\line(0, 1){20}}
\put(30, 70){\line(0, 1){20}}
\put(50, 50){\line(0, 1){20}}
\put(90, 50){\line(0, 1){40}}

\put(32, 78){$\scriptstyle{\bigl\}\om_i}$}
\put(52, 58){$\scriptstyle{\bigl\}\om_{i+1}}$}

\end{picture}
\begin{picture}(110, 120)
\put(30, 60){\makebox(50, 12){$\underrightarrow{r(\om_i,\,\om_{i+1})\ }$}}
\end{picture}
\begin{picture}(100, 120)(0, -10)

\linethickness{0.2pt}

\put(0, 10){\line(1, 0){100}}
\put(0, 30){\line(1, 0){100}}
\put(0, 50){\line(1, 0){100}}
\put(0, 70){\line(1, 0){30}}
\put(90, 70){\line(1, 0){10}}
\put(0, 90){\line(1, 0){10}}

\put(10, 0){\line(0, 1){90}}
\put(30, 0){\line(0, 1){70}}
\put(50, 0){\line(0, 1){50}}
\put(70, 0){\line(0, 1){50}}
\put(90, 0){\line(0, 1){50}}

\put(29, -7){$\scriptstyle{i}$}
\put(45, -7){$\scriptstyle{i\!+\!1}$}

\linethickness{2pt}

\put(0, 110){\line(1, 0){10}}
\put(10, 90){\line(1, 0){40}}
\put(50, 50){\line(1, 0){40}}
\put(90, 90){\line(1, 0){10}}

\put(10, 90){\line(0, 1){20}}
\put(50, 50){\line(0, 1){40}}
\put(90, 50){\line(0, 1){40}}

\put(30, 70){\dashbox{2}(20, 20){}}

\end{picture}
\end{center}
\caption{A possible move\label{fig:elso}}
\end{figure}\clearpage

\begin{table}[p]
\[
\begin{array}{|c||c|c||c|}
\hline
\text{with rate}&g_i\uparrow&h_i\uparrow&\text{a second class particle}\\
\hhline{|=#=|=#=|}
r(\ze_i,\,\ze_{i+1})-r(\eta_i,\,\ze_{i+1})&\bullet&&\curvearrowright\\
\hline
r(\eta_i,\,\eta_{i+1})-r(\eta_i,\,\ze_{i+1})&&\bullet&\curvearrowleft\\
\hline
r(\eta_i,\,\ze_{i+1})&\bullet&\bullet&\\
\hline
\end{array}
\]
\caption{Growth coupling rules}\label{tab:bas}
\end{table}\clearpage

\begin{table}[p]
\[
\begin{array}{|c||c|c|c|c||c|}
\hline
\text{with rate}&h_{i-1}\uparrow&g_{i-1}\uparrow&Q\curvearrowleft&S\curvearrowleft&\text{a second class particle}\\
\hhline{|=#=|=|=|=#=|}
f(-\om_i-1)-f(-\ze_i)&\bullet&&&&\curvearrowleft\\
\hline
\begin{array}{c}[f(-\om_i)-f(-\om_i-1)]-\\
-\frac{f(-\om_i)-f(-\ze_i)}{\ze_i-\om_i}\end{array}&\bullet&&\bullet&&\curvearrowleft\\
\hline
\frac{f(-\om_i)-f(-\ze_i)}{\ze_i-\om_i}&\bullet&&\bullet&\bullet&\curvearrowleft\\
\hline
f(-\ze_i)&\bullet&\bullet&&&\\
\hline
\end{array}
\]
\caption{Rates for $Q$ and $S$ to step left and for bricklayers at site $i=S=Q$ to lay brick on their left\label{tab:SQl}}
\end{table}\clearpage

\begin{table}[p]
\[
\begin{array}{|c||c|c|c|c||c|}
\hline
\text{with rate}&h_i\uparrow&g_i\uparrow&Q\curvearrowright&S\curvearrowright&\text{a second class particle}\\
\hhline{|=#=|=|=|=#=|}
\begin{array}{c}\frac{\ze_i-\om_i-1}{\ze_i-\om_i}\times\\
\times\left[f(\ze_i)-f(\om_i)\right]\end{array}&&\bullet&&&\curvearrowright\\
\hline
\begin{array}{c}\frac{f(\ze_i)-f(\om_i)}{\ze_i-\om_i}-\\
-\left[f(\om_i+1)-f(\om_i)\right]\end{array}&&\bullet&&\bullet&\curvearrowright\\
\hline
f(\om_i+1)-f(\om_i)&&\bullet&\bullet&\bullet&\curvearrowright\\
\hline
f(\om_i)&\bullet&\bullet&&&\\
\hline
\end{array}
\]
\caption{Rates for $Q$ and $S$ to step right and for bricklayers at site $i=S=Q$ to lay brick on their right\label{tab:SQr}}
\end{table}\clearpage

\begin{table}[p]
\[
\begin{array}{|c||c|c|c|c||c|}
\hline
\text{with rate}&h'_{i-1}\uparrow&g'_{i-1}\uparrow&Q\curvearrowleft&S'\curvearrowleft&\text{a second class particle}\\
\hhline{|=#=|=|=|=#=|}
\begin{array}{c}\frac{\ze'_i-\eta'_i-1}{\ze'_i-\eta'_i}\times\\
\times[f(-\eta'_i)-f(-\ze'_i)]\end{array}&\bullet&&&&\curvearrowleft\\
\hline
\begin{array}{c}\frac{f(-\eta'_i)-f(-\ze'_i)}{\ze'_i-\eta'_i}-\\
-[f(-\ze'_i+1)-f(-\ze'_i)]\end{array}&\bullet&&&\bullet&\curvearrowleft\\
\hline
f(-\ze'_i+1)-f(-\ze'_i)&\bullet&&\bullet&\bullet&\curvearrowleft\\
\hline
f(-\ze'_i)&\bullet&\bullet&&&\\
\hline
\end{array}
\]
\caption{Rates for $Q$ and $S'$ to step left and for bricklayers at site $i=S'=Q$ to lay brick on their left\label{tab:SQ'l}}
\end{table}\clearpage

\begin{table}[p]
\[
\begin{array}{|c||c|c|c|c||c|}
\hline
\text{with rate}&h'_i\uparrow&g'_i\uparrow&Q\curvearrowright&S'\curvearrowright&\text{a second class particle}\\
\hhline{|=#=|=|=|=#=|}
f(\ze'_i-1)-f(\eta'_i)&&\bullet&&&\curvearrowright\\
\hline
\begin{array}{c}\left[f(\ze'_i)-f(\ze'_i-1)\right]-\\
-\frac{f(\ze'_i)-f(\eta'_i)}{\ze'_i-\eta'_i}\end{array}&&\bullet&\bullet&&\curvearrowright\\
\hline
\frac{f(\ze'_i)-f(\eta'_i)}{\ze'_i-\eta'_i}&&\bullet&\bullet&\bullet&\curvearrowright\\
\hline
f(\eta'_i)&\bullet&\bullet&&&\\
\hline
\end{array}
\]
\caption{Rates for $Q$ and $S'$ to step right and for bricklayers at site $i=S'=Q$ to lay brick on their right\label{tab:SQ'r}}
\end{table}\clearpage
\end{document}